\documentclass[a4paper,10pt]{article}

\usepackage{amsthm,amsmath}
\usepackage{amssymb}
\usepackage{amsfonts}%,dsfont}
\usepackage[dvips]{epsfig}
\usepackage[dvips]{graphicx}
\usepackage{hyperref}
\usepackage[authoryear]{natbib}
\usepackage{hypernat}
\usepackage[T1]{fontenc}
\usepackage{shortcutssyl}
\usepackage{dsfont}

\theoremstyle{plain}

\theoremstyle{definition}

\newcommand{\egaldef}{:=} % egalite definissant la quantite de gauche
\newcommand{\flens}{\mapsto} % fleche d'application X->Y (ensembles)
\newcommand{\flapp}{\mapsto} % fleche d'application x->f(x) (elements)
\newcommand{\telque}{\, \mbox{ s.t. } \,} % tel que dans une definition d'ensemble

\newcommand{\set}[1]{\left\{ \left. #1 \right. \right\}}
\newcommand{\absj}[1]{\left\lvert #1 \right\rvert} %joli abs
\providecommand{\norm}[1]{\left \lVert #1 \right\rVert}

\newcommand{\guil}[1]{``#1''} %guillemets

\newcommand{\Prob}{\mathbb{P}} %probabilite
\newcommand{\Proba}{\mathbb{P}} %probabilite
\DeclareMathOperator{\var}{var} %variance
\newcommand{\sachant}{\, \right| \left. \,} % pour l'esp\'erance conditionnelle...
\newcommand{\bayes}{s}%estimateur Bayesien
\newcommand{\sha}{\widehat{s}^{\mathcal{A}}}%ERM
\newcommand{\mM}{m \in \M}

\newcommand{\mo}{\ensuremath{m^{\star}}}

\renewcommand{\sh}{\,\widehat{s}\,}
\renewcommand{\shm}{\,\widehat{s}_m\,}

\newcommand{\penid}{\pen_{\mathrm{id}}} % penalite ideale

\newcommand{\Los}[1]{\ensuremath{\mathcal{L}\paren{#1}}} % Risque
\newcommand{\LosP}[2]{\ensuremath{\mathcal{L}_{#2}\paren{#1}}} % Risque avec P en indice (2e argument)
\newcommand{\Loshnom}[1]{\ensuremath{\widehat{\mathcal{L}}^{\mathrm{#1}}}} % estimateur du risque dont le nom est precise en argument

\newcommand{\Loshval}{\Loshnom{H-O}} % estimateur du risque par validation
\newcommand{\Loshvc}{\Loshnom{CV}} % estimateur du risque par validation croisee
\newcommand{\Loshloo}{\Loshnom{LOO}} % estimateur du risque par leave-one-out
\newcommand{\Loshlpo}{\Loshnom{LPO}} % estimateur du risque par leave-p-out
\newcommand{\Loshvf}{\Loshnom{VF}} % estimateur du risque par validation croisee v-fold
\newcommand{\lh}{\ensuremath{\widehat{\lambda}}} % estimateur generique de lambda
\newcommand{\Il}{I_{\lambda}}
\newcommand{\lamm}{\lambda \in \Lambda_m}
\newcommand{\It}{\ensuremath{I^{(t)}}} % indices correspondant a l'echantillon d'entrainement (training)
\newcommand{\Iv}{\ensuremath{I^{(v)}}} % indices correspondant a l'echantillon de validation
\newcommand{\Dt}{\ensuremath{D_n^{(t)}}} % echantillon d'entrainement
\newcommand{\Dv}{\ensuremath{D_n^{(v)}}} % echantillon de validation
\renewcommand{\S}{\mathbb{S}} % ensemble des predicteurs / estimateurs / densites / ...

\newtheorem{postita}{Post-it}

\begin{document}

\title{A survey of cross-validation procedures for model selection}

\maketitle

\begin{center}
\author{
Sylvain Arlot,\\
CNRS ; Willow Project-Team, \\
Laboratoire d'Informatique de l'Ecole Normale Superieure\\
(CNRS/ENS/INRIA UMR 8548)\\
45, rue d'Ulm, 75\,230 Paris, France\\
\texttt{Sylvain.Arlot@ens.fr}\\
\ \\
Alain Celisse,\\
Laboratoire Paul Painlev\'e, UMR CNRS 8524,\\
Universit\'e des Sciences et Technologies de Lille 1\\
F-59\,655 Villeneuve d'Ascq Cedex, France\\
\texttt{Alain.Celisse@math.univ-lille1.fr}\\}
\end{center}

%%%%%%%%%%%%%%%%%%%%%%%%%%%%%%%%%%%%%%%%%%%%%%%%%%%%%%%%%%%%
%%%%%%%%%%%%%%%%%%%%%%%%%%%%%%%%%%%%%%%%%%%%%%%%%%%%%%%%%%%%

\begin{abstract}
Used to estimate the risk of an estimator or to perform model
selection, cross-validation is a widespread strategy because of its
simplicity and its apparent universality.
Many results exist on the model selection performances of
cross-validation procedures.
This survey intends to relate these results to the most recent
advances of model selection theory, with a particular emphasis on
distinguishing empirical statements from rigorous theoretical
results.
As a conclusion, guidelines are provided for choosing the best
cross-validation procedure according to the particular features of
the problem in hand.

\end{abstract}

%%%%%%%%%%%%%%%%%%%%%%%%%%%%%%%%%%%%%%%%%%%%%%%%%%%%%%%%%%%%
%%%%%%%%%%%%%%%%%%%%%%%%%%%%%%%%%%%%%%%%%%%%%%%%%%%%%%%%%%%%

\tableofcontents

%%%%%%%%%%%%%%%%%%%%%%%%%%%%%%%%%%%%%%%%%%%%%%%%%%%%%%%%%%%%%%%%%%%%%
%%%%%%%%%%%%%%%%%%%%%%%%%%%%%%%%%%%%%%%%%%%%%%%%%%%%%%%%%%%%%%%%%%%%%
%%%%%%%%%%%%%%%%%%%%%%%%%%%%%%%%%%%%%%%%%%%%%%%%%%%%%%%%%%%%%%%%%%%%%
%%%%%%%%%%%%%%%%%%%%%%%%%%%%%%%%%%%%%%%%%%%%%%%%%%%%%%%%%%%%%%%%%%%%%
%%%%%%%%%%%%%%%%%%%%%%% Introduction
%\input{intro.tex}% Introduction

\section{Introduction} \label{sec.intro}
Many statistical algorithms, such as likelihood maximization, least
squares and empirical contrast minimization, rely on the preliminary
choice of a model, that is of a set of parameters from which an
estimate will be returned. When several candidate models (thus
algorithms) are available, choosing one of them is called the model
selection problem.

Cross-validation (CV) is a popular strategy for model selection, and
more generally algorithm selection. The main idea behind CV is to
split the data (once or several times) for estimating the risk of
each algorithm: Part of the data (the training sample) is used for
training each algorithm, and the remaining part (the validation
sample) is used for estimating the risk of the algorithm. Then, CV
selects the algorithm with the smallest estimated risk.

Compared to the resubstitution error, CV avoids overfitting because
the training sample is independent from the validation sample (at
least when data are {\it i.i.d.}).
The popularity of CV mostly comes from the generality of the data
splitting heuristics, which only assumes that data are {\it i.i.d.}.
Nevertheless, theoretical and empirical studies of CV procedures do
not entirely confirm this ``universality''. Some CV procedures have
been proved to fail for some model selection problems, depending on
the goal of model selection: estimation or identification (see
Section~\ref{sec.modsel}). Furthermore, many theoretical questions
about CV remain widely open.

The aim of the present survey is to provide a clear picture of what
is known about CV, from both theoretical and empirical points of
view. More precisely, the aim is to answer the following questions:
What is CV doing? When does CV work for model selection, keeping in
mind that model selection can target different goals? Which CV
procedure should be used for each model selection problem?

The paper is organized as follows. First, the rest of
Section~\ref{sec.intro} presents the statistical framework. Although
non exhaustive, the present setting has been chosen general enough
for sketching the complexity of CV for model selection. The model
selection problem is introduced in Section~\ref{sec.modsel}. A brief
overview of some model selection procedures that are important to
keep in mind for understanding CV is given in
Section~\ref{sec.modselproc}.
The most classical CV procedures are defined in
Section~\ref{sec.def}.
Since they are the keystone of the behaviour of CV for model
selection, the main properties of CV estimators of the risk for a
fixed model are detailed in Section~\ref{sec.riskestim}. Then, the
general performances of CV for model selection are described, when
the goal is either estimation (Section~\ref{sec.cvefficient}) or
identification (Section~\ref{sec.cvconsistent}). Specific properties
of CV in some particular frameworks are discussed in
Section~\ref{sec.specific}.
Finally, Section~\ref{sec.complex} focuses on the algorithmic
complexity of CV procedures, and Section~\ref{sec.conclu} concludes
the survey by tackling several practical questions about CV.

\subsection{Statistical framework} \label{sec.def.cadre.pb}
Assume that some data $\xi_1, \ldots, \xi_n \in \Xi$ with common
distribution $P$ are observed. Throughout the paper---except in
Section~\ref{sec.dependent}---the $\xi_i$ are assumed to be
independent.
The purpose of statistical inference is to estimate from the data
$\paren{\xi_i}_{1 \leq i \leq n}$ some target feature $\bayes$ of the
unknown distribution $P$, such as the mean or the variance of $P$.
Let $\S$ denote the set of possible values for $s$.
The quality of $t \in \S$, as an approximation of $\bayes$, is
measured by its loss $\Los{t}$, where $\mathcal{L}: \S \flens \R$ is
called the {\em loss function}, and is assumed to be minimal for
$t=\bayes$. Many loss functions can be chosen for a given statistical
problem.

Several classical loss functions are defined by
\begin{equation}
\label{def.loss} \Los{t} = \LosP{t}{P} \egaldef \E_{\xi \sim P}
\croch{ \gamma\paren{ t ; \xi }} \enspace ,
\end{equation}
where $\gamma: \S \times \Xi \flens [0,\infty)$ is called a {\em
contrast function}.
Basically, for $t \in \S$ and $\xi \in \Xi$, $\gamma(t;\xi)$ measures
how well $t$ is in accordance with observation of $\xi$, so that the
loss of $t$, defined by \eqref{def.loss}, measures the average
accordance between $t$ and new observations $\xi$ with distribution
$P$.
Therefore, several frameworks such as transductive learning do not
fit definition \eqref{def.loss}.
Nevertheless, as detailed in Section~\ref{sec.def.cadre.ex},
definition \eqref{def.loss} includes most classical statistical
frameworks.

Another useful quantity is the {\em excess loss}
\[ \perte{t} \egaldef \LosP{t}{P} - \LosP{\bayes}{P} \geq 0\enspace, \]
which is related to the risk of an estimator $\sh$ of the target $s$
by
\begin{align*}
R(\sh)= \E_{\xi_1,\ldots,\xi_n\sim P}\croch{\perte{\sh}}\enspace .
\end{align*}

\subsection{Examples} \label{sec.def.cadre.ex}
The purpose of this subsection is to show that the framework of
Section~\ref{sec.def.cadre.pb} includes several important statistical
frameworks. This list of examples does not pretend to be exhaustive.
\paragraph{Density estimation} aims at estimating the density $\bayes$
of $P$ with respect to some given measure $\mu$ on $\Xi$. Then, $\S$
is the set of densities on $\Xi$ with respect to $\mu$.
For instance, taking $\gamma(t;x) = - \ln(t(x))$ in \eqref{def.loss},
the loss is minimal when $t = \bayes$ and the excess loss
\[ \perte{t} = \LosP{t}{P} - \LosP{\bayes}{P} = \E_{\xi \sim P} \croch{ \ln \paren{ \frac{ \bayes(\xi)} {t(\xi)} } } = \int \bayes \ln \paren{ \frac{s}{t} } \textup{d}\mu \]
is the Kullback-Leibler divergence between distributions $t \mu$ and
$\bayes \mu$.

\paragraph{Prediction} aims at predicting a quantity of interest $Y \in \Y$ given
an explanatory variable $X \in \X$ and a sample of observations
$(X_1,Y_1), \ldots, (X_n, Y_n)$. In other words, $\Xi = \X \times
\Y$, $\S$ is the set of measurable mappings $\X \flens \Y$ and the
contrast $\gamma(t;(x,y))$ measures the discrepancy between the
observed $y$ and its predicted value $t(x)$.
Two classical prediction frameworks are regression and
classification, which are detailed below.

\paragraph{Regression} corresponds to continuous $\Y$, that is $\Y \subset \R$ (or $\R^k$ for multivariate regression), the feature space $\X$ being typically a subset of $\R^{\ell}$.
Let $\bayes$ denote the regression function, that is $\bayes(x) =
\E_{(X,Y) \sim P} \croch{ Y \sachant X = x}$, so that
\[ \forall i, \quad Y_i = \bayes(X_i) + \epsilon_i \qquad \mbox{with} \quad  \E\croch{\epsilon_i \sachant X_i} = 0 \enspace . \]
A popular contrast in regression is the {\em least-squares contrast}
$ \gamma\paren{t;(x,y)} = (t(x) - y)^2 $, which is minimal over $\S$
for $t = \bayes$, and the excess loss is
\[ \perte{t} = \E_{(X,Y) \sim P} \croch{ \paren{ s(X) - t(X) }^2 } \enspace . \]
Note that the excess loss of $t$ is the square of the $L^2$ distance
between $t$ and $s$, so that prediction and  estimation are
equivalent goals.

\paragraph{Classification} corresponds to finite $\Y$ (at least discrete). In particular, when $\Y = \set{0,1}$, the prediction problem is called {\em binary (supervised) classification}.
With the 0-1 contrast function $\gamma(t;(x,y)) = \1_{t(x) \neq y}$,
the minimizer of the loss is the so-called Bayes classifier $\bayes$
defined by
\[ s(x) = \1_{\eta(x) \geq 1/2} \enspace , \]
where $\eta$ denotes the regression function $\eta(x) = \Proba_{(X,Y)
\sim P}\paren{ Y=1 \sachant X = x}$.

Remark that a slightly different framework is often considered in
binary classification.
Instead of looking only for a classifier, the goal is to estimate
also the  confidence in the classification made at each point: $\S$
is the set of measurable mappings $\X \flens \R$, the classifier $x
\flapp \1_{t(x) \geq 0}$ being associated to any $t \in \S$.
Basically, the larger $\absj{t(x)}$, the more confident we are in the
classification made from $t(x)$.
A classical family of losses associated with this problem is defined
by \eqref{def.loss} with the contrast $\gamma_{\phi} \paren{t; (x,y)}
= \phi \paren{ - (2y-1) t(x) }$ where $\phi : \R \flens [0,  \infty)$
is some function. The 0-1 contrast corresponds to $\phi(u) = \1_{u
\geq 0}$. The convex loss functions correspond to the case where
$\phi$ is convex, nondecreasing with $\lim_{- \infty} \phi = 0$ and
$\phi(0)=1$.
Classical examples are $\phi(u) = \max\set{ 1 + u, 0}$ (hinge),
$\phi(u) = \exp(u)$, and $\phi(u) = \log_2 \paren{ 1 + \exp(u)}$
(logit). The corresponding losses are used as objective functions by
several classical learning algorithms such as support vector machines
(hinge) and boosting (exponential and logit).

Many references on classification theory, including model selection,
can be found in the survey by \cite{Bou_Bou_Lug:2005}.

\subsection{Statistical algorithms}  \label{sec.def.cadre.sol}
\sloppy In this survey, a {\em statistical algorithm} $\A$ is any
(measurable) mapping $\A: \bigcup_{n \in \N} \Xi^n \flens S$. The
idea is that data $D_n = \paren{\xi_i}_{1 \leq i \leq n} \in \Xi^n$
will be used as an input of $\A$, and that the output of $\A$,
$\A(D_n)=\sha(D_n) \in \S$, is an estimator of $s$.
The quality of $\A$ is then measured by $\LosP{\sha(D_n)}{P}$, which
should be as small as possible.
In the sequel, the algorithm $\A$ and the estimator $\sha(D_n)$ are
often identified when no confusion is possible.

\medskip

{\em Minimum contrast estimators} form a classical family of
statistical algorithms, defined as follows. Given some subset $S$ of
$\S$ that we call a {\em model}, a minimum contrast estimator of $s$
is any minimizer of the empirical contrast
\[ t \flapp \LosP{t}{P_n} = \frac{1}{n} \sum_{i=1}^n \gamma\paren{t ; \xi_i}, \qquad \mbox{where} \quad P_n = \frac{1}{n} \sum_{i=1}^n \delta_{\xi_i}\enspace, \] over $S$.
The idea is that the empirical contrast $\LosP{t}{P_n}$ has an
expectation $\LosP{t}{P}$ which is minimal over $\S$ at $s$. Hence,
minimizing $\LosP{t}{P_n}$ over a set $S$ of candidate values for $s$
hopefully leads to a good estimator of $s$.
Let us now give three popular examples of empirical contrast
minimizers:
\begin{itemize}
\item {\em Maximum likelihood estimators}: take $\gamma(t;x) = - \ln(t(x))$ in the density estimation setting.
A classical choice for $S$ is the set of piecewise constant functions
on a regular partition of $\Xi$ with $K$ pieces.
\item {\em Least-squares estimators}: take $\gamma(t;(x,y)) = (t(x) - y)^2$ the least-squares contrast in the regression setting.
For instance, $S$ can be the set of piecewise constant functions on
some fixed partition of $\X$ (leading to regressograms), or a vector
space spanned by the first vectors of wavelets or Fourier basis,
among many others.
Note that regularized least-squares algorithms such as the Lasso,
ridge regression and spline smoothing also are least-squares
estimators, the model $S$ being some ball of a (data-dependent)
radius for the $L^1$ (resp. $L^2$) norm in some high-dimensional
space. Hence, tuning the regularization parameter for the LASSO or
SVM, for instance, amounts to perform model selection from a
collection of models.
\item {\em Empirical risk minimizers}, following the terminology of \cite{Vap:1982}: take any contrast function $\gamma$ in the prediction setting.
When $\gamma$ is the 0-1 contrast, popular choices for $S$ lead to
linear classifiers, partitioning rules, and neural networks.
Boosting and Support Vector Machines classifiers also are empirical
contrast minimizers over some data-dependent model $S$, with contrast
$\gamma = \gamma_{\phi}$ for some convex functions $\phi$.
\end{itemize}

\medskip

Let us finally mention that many other classical statistical
algorithms can be considered with CV, for instance
local average estimators in the prediction framework such as
$k$-Nearest Neighbours and Nadaraya-Watson kernel estimators.
The focus will be mainly kept on minimum contrast estimators to keep
the length of the survey reasonable.

%%%%%%%%%%%%%%%%%%%%%%%%%%%%%%%%%%%%%%%%%%%%%%%%%%%%%%%%%%%%%%%%%%%%%%%%%%%%
%%%%%%%%%%%%%%%%%%%%%%%%%%%%%%%%%%%%%%%%%%%%%%%%%%%%%%%%%%%%%%%%%%%%%%%%%%%%
%%%%%%%%%%%%%%%%%%%%%%%%%%%%%%%%%%%%%%%%%%%%%%%%%%%%%%%%%%%%%%%%%%%%%%%%%%%%
%%%%%%%%%%%%%%%%%%%%%%%%%%%%%%%%%%%%%%%%%%%%%%%%%%%%%%%%%%%%%%%%%%%%%%%%%%%%
%%%%%%%%%%%%%%%%%%%%%%%%%%%%%%%%%%%%%%%%%%%%%%%%%%%%%%%%%%%%%%%%%%%%%%%%%%%%

%%%%%%%%%%%%%%%%%%%%%%%%%%%%%%%%%%%%%%%%%%%%%%%%%%%%%%%%%%%%%%%%%%%%%%%%%%%%
%%%%%%%%%%%%%%%%%%%%%%% Model selection
%\input{modsel.tex}% Model selection

\section{Model selection} \label{sec.modsel}
Usually, several statistical algorithms can be used for solving a
given statistical problem. Let $\paren{\ERM_{\lambda}}_{\lambda \in
\Lambda}$ denote such a family of candidate statistical algorithms.
The {\em algorithm selection problem} aims at choosing from data one
of these algorithms, that is, choosing some $\lh(D_n) \in \Lambda$.
Then, the final estimator of $s$ is given by $\ERM_{\lh(D_n)} (D_n)$.
The main difficulty is that the same data are used for training the
algorithms, that is, for computing $\paren{\ERM_{\lambda} (D_n)
}_{\lambda \in \Lambda}$, and for choosing $\lh(D_n)$ .

\subsection{The model selection paradigm} \label{sec.modsel.paradigm}
Following Section~\ref{sec.def.cadre.sol}, let us focus on the {\em
model selection problem}, where candidate algorithms are minimum
contrast estimators and the goal is to choose a model $S$.
Let $\paren{S_m}_{\mM_n}$ be a family of models, that is, $S_m
\subset \S$. Let $\gamma$ be a fixed contrast function, and for every
$\mM_n$, let $\ERM_m$ be a minimum contrast estimator over model
$S_m$ with contrast $\gamma$. The goal is to choose $\mh(D_n) \in
\M_n$ from data only.

\medskip

The choice of a model $S_m$ has to be done carefully. Indeed, when
$S_m$ is a \guil{small} model, $\ERM_m$ is a poor statistical
algorithm except when $s$ is very close to $S_m$, since
\[ \perte{\ERM_m} \geq \inf_{t \in S_m} \set{ \perte{t} } \egaldef \perte{S_m} \enspace .\]
The lower bound $\perte{S_m}$ is called the {\em bias} of model
$S_m$, or {\em approximation error}. The bias is a nonincreasing
function of $S_m$.

On the contrary, when $S_m$ is \guil{huge}, its bias $\perte{S_m}$ is
small for most targets $s$, but $\ERM_m$ clearly overfits. Think for
instance of $S_m$ as the set of all continuous functions on $[0,1]$
in the regression framework.
More generally, if $S_m$ is a vector space of dimension $D_m$, in
several classical frameworks,
\begin{equation} \label{eq.biais-var} \E\croch{ \perte{\ERM_m(D_n)} } \approx \perte{S_m} + \lambda D_m \end{equation}
where $\lambda>0$ does not depend on $m$.
For instance, $\lambda = 1 / (2n)$ in density estimation using the
likelihood contrast, and $\lambda = \sigma^2 /n$ in regression using
the least-squares contrast and assuming $\var \paren{ Y\sachant X } =
\sigma^2$ does not depend on $X$.
The meaning of \eqref{eq.biais-var} is that a good model choice
should balance the bias term $\perte{S_m}$ and the {\em variance}
term $\lambda D_m$, that is solve the so-called {\em bias-variance
trade-off}.
By extension, the variance term, also called {\em estimation error},
can be defined by
\[ \E\croch{ \perte{\ERM_m(D_n)} } - \perte{S_m} =
\E\croch{\LosP{\shm}{P}} - \inf_{t\in S_m}\LosP{t}{P}\enspace,\] even
when \eqref{eq.biais-var} does not hold.

The interested reader can find a much deeper insight into model
selection in the Saint-Flour lecture notes by
\cite{Mas:2003:St-Flour}.

\medskip

Before giving examples of classical model selection procedures, let
us mention the two main different goals that model selection can
target: estimation and identification.

\subsection{Model selection for estimation} \label{sec.modsel.estim}
On the one hand, the goal of model selection is {\em estimation} when
$\ERM_{\mh(D_n)} (D_n)$ is used as an approximation of the target
$s$, and the goal is to minimize its loss. For instance, AIC and
Mallows' $C_p$ model selection procedures are built for estimation
(see Section~\ref{sec.modselproc.unbiased}).

The quality of a model selection procedure $D_n \flapp \mh(D_n)$,
designed for estimation, is measured by the excess loss of
$\ERM_{\mh(D_n)} (D_n)$. Hence, the best possible model choice for
estimation is the so-called {\em oracle} model $S_{\mo}$, defined by
\begin{equation} \label{eq.mo}
\mo = \mo(D_n) \in \arg\min_{\mM_n} \set{ \perte{ \ERM_m (D_n) } }
\enspace .
\end{equation}
Since $\mo(D_n)$ depends on the unknown distribution $P$ of data, one
cannot expect to select $\mh(D_n)=\mo(D_n)$ almost surely.
Nevertheless, we can hope to select $\mh(D_n)$ such that
$\sh_{\mh(D_n)}$ is almost as close to $s$ as $\sh_{\mo(D_n)}$.
Note that there is no requirement for $s$ to belong to
$\bigcup_{m\in\mathcal{M}_n}S_m$.

\medskip

Depending on the framework, the optimality of a model selection
procedure for estimation is assessed in at least two different ways.

First, in the asymptotic framework, a model selection procedure $\mh$
is called {\em efficient} (or asymptotically optimal)
 when it leads to $\mh$ such that
\[ \frac{ \perte{ \ERM_{\mh(D_n)} (D_n) } } { \inf_{\mM_n} \set{ \perte{ \ERM_m (D_n) }} } \xrightarrow[n \rightarrow \infty]{a.s.} 1 \enspace . \]
Sometimes, a weaker result is proved, the convergence holding only in
probability.

Second, in the non-asymptotic framework, a model selection procedure
satisfies an {\em oracle inequality} with constant $C_n \geq 1$ and
remainder term $R_n \geq 0$ when
\begin{equation} \label{eq.oracle-ineq}
\perte{ \ERM_{\mh(D_n)} (D_n) } \leq C_n \inf_{\mM_n} \set{ \perte{
\ERM_m (D_n) } } + R_n
\end{equation}
holds either in expectation or with large probability (that is, a
probability larger than $1 - C^{\prime}/ n^2$, for some positive
constant $C^{\prime}$).
Note that if \eqref{eq.oracle-ineq} holds on a large probability
event with $C_n$ tending to~1 when $n$ tends to infinity and $R_n \ll
\perte{\ERM_{\mo}(D_n)}$, then the model selection procedure $\mh$ is
efficient.

\medskip

In the estimation setting, model selection is often used for building {\em adaptive estimators}, assuming that $s$ belongs to some function space $\mathcal{T}_{\alpha}$ \citep{Bar_Bir_Mas:1999}.%\citep{Bir_Mas:1997}.
Then, a model selection procedure $\mh$ is optimal when it leads to
an estimator $\ERM_{\mh(D_n)} (D_n)$ (approximately) minimax with
respect to $\mathcal{T}_{\alpha}$ without knowing $\alpha$, provided
the family $\paren{S_m}_{\mM_n}$ has been well-chosen.

\subsection{Model selection for identification} \label{sec.modsel.identif}
On the other hand, model selection can aim at identifying the ``true
model'' $S_{m_0}$, defined as the ``smallest'' model among
$\paren{S_m}_{\mM_n}$ to which $s$ belongs.
In particular, $s \in \bigcup_{\mM_n} S_m$ is assumed in this
setting.
A typical example of model selection procedure built for
identification is BIC (see Section~\ref{sec.modselproc.identif}).

The quality of a model selection procedure designed for
identification is measured by its probability of recovering the true
model $m_0$.
Then, a model selection procedure is called {\em (model) consistent}
when
\[ \Proba \paren{ \mh(D_n) = m_0 } \xrightarrow[n \rightarrow \infty]{} 1 \enspace .\]
Note that identification can naturally be extended to the general
algorithm selection problem, the ``true model'' being replaced by the
statistical algorithm whose risk converges at the fastest rate
\citep[see for instance][]{Yan:2007b}.

\subsection{Estimation {\em vs.} identification} \label{sec.modsel.estim-vs-identif}
When a true model exists, model consistency is clearly a stronger
property than efficiency defined in Section~\ref{sec.modsel.estim}.
However, in many frameworks, no true model does exist so that
efficiency is the only well-defined property.

Could a model selection procedure be model consistent in the former
case (like BIC) and efficient in the latter case (like AIC)?
The general answer to this question, often called the AIC-BIC
dilemma, is negative: \cite{Yan:2005a} proved in the regression
framework that no model selection procedure can be simultaneously
model consistent and minimax rate optimal.
Nevertheless, the strengths of AIC and BIC can sometimes be shared;
see for instance the introduction of a paper by \cite{Yan:2005a} and
a recent paper by \cite{Erv_Gru_Roo:2008}.

\section{Overview of some model selection procedures} \label{sec.modselproc} %\label{sec.modsel.other}
Several approaches can be used for model selection. Let us briefly
sketch here some of them, which are particularly helpful for
understanding how CV works. Like CV, all the procedures considered in
this section select
\begin{equation} \label{eq.uerp}
\mh(D_n) \in \arg\min_{\mM_n} \set{ \crit(m; D_n)} \enspace ,
\end{equation}
where $\forall \mM_n$, $\crit(m; D_n)=\crit(m) \in \R$ is some
data-dependent criterion.

A particular case of \eqref{eq.uerp} is {\em penalization}, which
consists in choosing the model minimizing the sum of empirical
contrast and some measure of complexity of the model (called penalty)
which can depend on the data, that is,
\begin{equation} \label{eq.pen}
\mh (D_n) \in \arg\min_{\mM_n} \set{ \LosP{\ERM_{m}}{P_n} + \pen(m;
D_n) } \enspace .
\end{equation}
This section does not pretend to be exhaustive. Completely different
approaches exist for model selection, such as the Minimum Description
Length (MDL) \citep{Riss83}, and the Bayesian approaches. The
interested reader will find more details and references on model
selection procedures in the books by \cite{Bur_And:2002} or
\cite{Mas:2003:St-Flour} for instance.

Let us focus here on five main categories of model selection
procedures, the first three ones coming from a classification made by
\cite{Sha:1997} in the linear regression framework.

%%%%%%%%%%%%%%%%%%%
%%% 1ere categorie
%%%%%%%%%%%%%%%%%%%
%%%
\subsection{The unbiased risk estimation principle} \label{sec.modselproc.unbiased}
When the goal of model selection is estimation, many model selection
procedures are of the form \eqref{eq.uerp} where $\crit(m;D_n)$
unbiasedly estimates (at least, asymptotically) the loss
$\LosP{\ERM_m}{P}$.
This general idea is often called unbiased risk estimation principle,
or Mallows' or Akaike's heuristics.

In order to explain why this strategy can perform well, let us write
the starting point of most theoretical analysis of procedures defined
by \eqref{eq.uerp}: By definition \eqref{eq.uerp}, for every $\mM_n$,
\begin{equation} \label{eq.start-oracle}
\perte{\ERM_{\mh}} + \crit(\mh) - \LosP{\ERM_{\mh}}{P} \leq
\perte{\ERM_{m}} + \crit(m) - \LosP{\ERM_{m}}{P} \enspace .
\end{equation}
If $\E\croch{ \crit(m) - \LosP{\ERM_{m}}{P} } = 0$ for every $\mM_n$,
then concentration inequalities are likely to prove that
$\varepsilon_n^-, \varepsilon_n^+ > 0$ exist such that
\[ \forall \mM_n, \quad  \varepsilon_n^+ \geq \frac{ \crit(m) - \LosP{\ERM_{m}}{P} }{ \perte{\ERM_m} } \geq  - \varepsilon_n^- > -1 \]
with high probability, at least when $\card(\M_n) \leq C n^{\alpha}$
for some $C,\alpha \geq 0$.
Then, \eqref{eq.start-oracle} directly implies an oracle inequality
like \eqref{eq.oracle-ineq} with $C_n = (1 + \varepsilon_n^+) / (1 -
\varepsilon_n^-)$.
If $\varepsilon_n^+, \varepsilon_n^- \rightarrow 0$ when $n
\rightarrow \infty$, this proves the procedure defined by
\eqref{eq.uerp} is efficient.

\medskip

Examples of model selection procedures following the unbiased risk
estimation principle are FPE \citep[Final Prediction
Error,][]{Aka:1969}, several cross-validation procedures including
the Leave-one-out (see Section~\ref{sec.def}), and GCV
\citep[Generalized Cross-Validation,][see
Section~\ref{sec.def.classex.other}]{Cra_Wah:1979}.
%%%%%
%
With the penalization approach \eqref{eq.pen}, the unbiased risk
estimation principle is that $\E\croch{\pen(m)}$ should be close to
the ``ideal penalty''
\[ \penid(m) \egaldef \LosP{\ERM_m}{P} - \LosP{\ERM_m}{P_n} \enspace . \]
Several classical penalization procedures follow this principle, for
instance:
\begin{itemize}
\item With the log-likelihood contrast, AIC \citep[Akaike Information
Criterion,][]{Aka:1973} and its corrected versions
\citep{Sug:1978,Hur_Tsa:1989}.
\item With the least-squares contrast, Mallows' $C_p$ \citep{Mal:1973} and several refined versions of $C_p$ \citep[see for instance][]{Bar:2002}.
\item With a general contrast, covariance penalties \citep{Efr:2004}.
\end{itemize}

AIC, Mallows' $C_p$ and related procedures have been proved to be
optimal for estimation in several frameworks, provided $\card(\M_n)
\leq C n^{\alpha}$ for some constants $C,\alpha \geq 0$ \citep[see
the paper by][and references therein]{Bir_Mas:2006}.

\medskip

%%%%%%%%%%%%%% Resampling

The main drawback of penalties such as AIC or Mallows' $C_p$ is their
dependence on some assumptions on the distribution of data. For
instance, Mallows' $C_p$ assumes the variance of $Y$ does not depend
on $X$. Otherwise, it has a suboptimal performance \citep{Arlo08}.

Several resampling-based penalties have been proposed to overcome
this problem, at the price of a larger computational complexity, and
possibly slightly worse performance in simpler frameworks; see a
paper by \cite{Efr:1983} for bootstrap, and a paper by \cite{Arlo08c}
and references therein for generalization to exchangeable weights.

\medskip

\sloppy Finally, note that all these penalties depend on multiplying
factors which are not always known (for instance, the noise-level,
for Mallows' $C_p$).
\cite{Bir_Mas:2006} proposed a general data-driven procedure for
estimating such multiplying factors, which satisfies an oracle
inequality with $C_n \rightarrow 1$ in regression \citep[see
also][]{Arl_Mas:2009:pente}.

%%%%%%%%%%%%%%%%%%%
%%% 2eme categorie
%%%%%%%%%%%%%%%%%%%
%%%
\subsection{Biased estimation of the risk} \label{sec.modselproc.biased}
Several model selection procedures are of the form \eqref{eq.uerp}
where $\crit(m)$ does not unbiasedly estimate the loss
$\LosP{\ERM_m}{P}$: The weight of the variance term compared to the
bias in $\E\croch{\crit(m)}$ is slightly larger than in the
decomposition \eqref{eq.biais-var} of $\LosP{\ERM_m}{P}$.
From the penalization point of view, such procedures are {\em
overpenalizing}.

Examples of such procedures are FPE$_{\alpha}$ \citep{Bha_Dow:1977}
and GIC$_{\lambda}$ \citep[Generalized Information
Criterion,][]{Nis:1984,Sha:1997} with $\alpha,\lambda>2$, which are
closely related. Some cross-validation procedures, such as
Leave-$p$-out with $p/n \in (0,1)$ fixed, also belong to this
category (see Section~\ref{subsec.classical.examples.exhaust}).
Note that FPE$_{\alpha}$ with $\alpha=2$ is FPE, and GIC$_{\lambda}$
with $\lambda=2$ is close to FPE and Mallows' $C_p$.

\medskip

When the goal is estimation, there are two main reasons for using
``biased'' model selection procedures.
First, experimental evidence show that overpenalizing often yields
better performance when the signal-to-noise ratio is small \citep[see
for instance][Chapter~11]{Arlo07}.

Second, when the number of models $\card(\M_n)$ grows faster than any
power of $n$, as in the complete variable selection problem with $n$
variables, then the unbiased risk estimation principle fails. From
the penalization point of view, \cite{Bir_Mas:2006} proved that when
$\card(\M_n) = e^{\kappa n}$ for some $\kappa >0$, the minimal amount
of penalty required so that an oracle inequality holds with $C_n =
\grandO(1)$ is much larger than $\penid(m)$.
In addition to the FPE$_{\alpha}$ and GIC$_{\lambda}$ with suitably
chosen $\alpha, \lambda$, several penalization procedures have been
proposed for taking into account the size of $\M_n$
\citep{Bar_Bir_Mas:1999,Bar:2002,Bir_Mas:2002,Sau:2006}. In the same
papers, these procedures are proved to satisfy oracle inequalities
with $C_n$ as small as possible, typically of order $\ln(n)$ when
$\card(\M_n) = e^{\kappa n}$.
%
%%%%% Pourquoi surpenaliser? SNR petit
%

%%%%%%%%%%%%%%%%%%%
%%% 3eme categorie
%%%%%%%%%%%%%%%%%%%
%%%
\subsection{Procedures built for identification} \label{sec.modselproc.identif}
Some specific model selection procedures are used for identification.
A typical example is BIC \citep[Bayesian Information
Criterion,][]{Sch:1978}.

More generally, \cite{Sha:1997} showed that several procedures
identify consistently the correct model in the linear regression
framework as soon as they overpenalize within a factor tending to
infinity with $n$, for instance, GIC$_{\lambda_n}$ with $\lambda_n
\rightarrow + \infty$, FPE$_{\alpha_n}$ with $\alpha_n \rightarrow
+\infty$ \citep{Shi:1984}, and several CV procedures such as
Leave-$p$-out with $p=p_n\sim n$.
BIC is also part of this picture, since it coincides with
GIC$_{\ln(n)}$.

In another paper, \cite{Sha:1996} showed that $m_n$-out-of-$n$
bootstrap penalization is also model consistent as soon as $m_n \sim
n$. Compared to Efron's bootstrap penalties, the idea is to estimate
$\penid$ with the $m_n$-out-of-$n$ bootstrap instead of the usual
bootstrap, which results in overpenalization within a factor tending
to infinity with $n$ \citep{Arlo08c}.

Most MDL-based procedures can also be put into this category of model
selection procedures \citep[see][]{Gru:2007}.
Let us finally mention the Lasso \citep{Tibs96}
%, its variants \citep[e.g.][]{Bac:2008:GL-MKL,Bac:2009}
and other $\ell^1$ penalization procedures, which have recently
attracted much attention \citep[see for instance][]{Hes_etal:2008}.
They are a computationally efficient way of identifying the true
model in the context of variable selection with many variables.

%%%%%%%%%%%%%%%%%%%
%%% 4eme categorie
%%%%%%%%%%%%%%%%%%%
%%%
\subsection{Structural risk minimization} \label{sec.modselproc.Vapnik}
In the context of statistical learning, %in particular classification,
\cite{Vap_Cer:1974} proposed the structural risk minimization
approach \citep[see also][]{Vap:1982,Vap:1998}.
Roughly, the idea is to penalize the empirical contrast with a
penalty (over)-estimating
\[ \penidglo(m) \egaldef
\sup_{t \in S_m} \set{ \LosP{t}{P} - \LosP{t}{P_n} } \geq \penid(m)
\enspace . \]
Such penalties have been built using the Vapnik-Chervonenkis
dimension, the combinatorial entropy, (global) Rademacher
complexities \citep{Kol:2001,Bar_Bou_Lug:2002}, (global) bootstrap
penalties \citep{Fro:2007}, Gaussian complexities or the maximal
discrepancy \citep{Bar_Men:2002}. These penalties are often called
{\em global} because $\penidglo(m)$ is a supremum over $S_m$.

The localization approach \citep[see][]{Bou_Bou_Lug:2005} has been
introduced in order to obtain penalties closer to $\penid$ (such as
local Rademacher complexities), hence smaller prediction errors when
possible \citep{Bar_Bou_Men:2005,Kol:2006}. Nevertheless, these
penalties are still larger than $\penid(m)$ and can be difficult to
compute in practice because of several unknown constants.

A non-asymptotic analysis of several global and local penalties can
be found in the book by \cite{Mas:2003:St-Flour} for instance; see
also \cite{Kol:2006} for recent results on local penalties.

\subsection{{\em Ad hoc} penalization} \label{sec.modselproc.adhoc}
Let us finally mention that penalties can also be built according to
particular features of the problem.
For instance, penalties can be proportional to the $\ell^p$ norm of
$\ERM_m$ (similarly to $\ell^p$-regularized learning algorithms) when
having an estimator with a controlled $\ell^p$ norm seems better.
The penalty can also be proportional to the squared norm of $\ERM_m$
in some reproducing kernel Hilbert space (similarly to kernel ridge
regression or spline smoothing), with a kernel adapted to the
specific framework.
More generally, any penalty can be used, as soon as $\pen(m)$ is
larger than the estimation error (to avoid overfitting) and the best
model for the final user is not the oracle $\mo$, but more like
\[ \arg\min_{\mM_n} \set{ \perte{ S_m  } + \kappa \pen(m) } \]
for some $\kappa >0$.

\subsection{Where are cross-validation procedures in this picture?} \label{sec.modselproc.CV}
The family of CV procedures, which will be described and deeply
investigated in the next sections, contains procedures in the first
three categories.
CV procedures are all of the form \eqref{eq.uerp}, where $\crit(m)$
either estimates (almost) unbiasedly the loss $\LosP{\ERM_m}{P}$, or
overestimates the variance term (see
Section~\ref{sec.modsel.paradigm}). In the latter case, CV procedures
either belong to the second or the third category, depending on the
overestimation level.

This fact has two major implications.
First, CV itself does not take into account prior information for
selecting a model. To do so, one can either add to the CV estimate of
the risk a penalty term (such as $\norm{\ERM_m}_p$), or use prior
information to pre-select a subset of models $\Mt(D_n) \subset \M_n$
before letting CV select a model among $\paren{S_m}_{m \in
\Mt(D_n)}$.

Second, in statistical learning, CV and resampling-based procedures
are the most widely used model selection procedures. Structural risk
minimization is often too pessimistic, and other alternatives rely on
unrealistic assumptions.
But if CV and resampling-based procedures are the most likely to
yield good prediction performances, their theoretical grounds are not
that firm, and too few CV users are careful enough when choosing a CV
procedure to perform model selection.
Among the aims of this survey is to point out both positive and
negative results about the model selection performance of CV.

%%%%%%%%%%%%%%%%%%%%%%%%%%%%%%%%%%%%%%%%%%%%%%%%%%%%%%%%%%%%%%%%%%%%%%%%%
%%%%%%%%%%%%%%%%%%%%%%%%%%%%%%%%%%%%%%%%%%%%%%%%%%%%%%%%%%%%%%%%%%%%%%%%%
%%%%%%%%%%%%%%%%%%%%%%%%%%%%%%%%%%%%%%%%%%%%%%%%%%%%%%%%%%%%%%%%%%%%%%%%%
%%%%%%%%%%%%%%%%%%%%%%%%%%%%%%%%%%%%%%%%%%%%%%%%%%%%%%%%%%%%%%%%%%%%%%%%%
%%%%%%%%%%%%%%%%%%%%% Cross-validation
%\input{crossval.tex}% Cross-validation procedures

\section{Cross-validation procedures} \label{sec.def}
The purpose of this section is to describe the rationale behind CV
and to define the different CV procedures.
Since all CV procedures are of the form \eqref{eq.uerp}, defining a
CV procedure amounts to define the corresponding CV estimator of the
risk of an algorithm $\A$, which will be $\crit(\cdot)$ in
\eqref{eq.uerp}.

\subsection{Cross-validation philosophy} \label{sec.def.philo}
As noticed in the early 30s by \cite{Lars31}, training an algorithm
and evaluating its statistical performance on the same data yields an
overoptimistic result. CV was raised to fix this issue
\citep{MoTu68,Sto:1974,Gei:1975}, starting from the remark that
testing the output of the algorithm on new data would yield a good
estimate of its performance \citep{Brei96c}.

In most real applications, only a limited amount of data is
available, which led to the idea of {\em splitting the data}: Part of
the data (the training sample) is used for training the algorithm,
and the remaining data (the validation sample) is used for evaluating
its performance. The validation sample can play the role of new data
as soon as data are {\it i.i.d.}.

Data splitting yields the {\em validation} estimate of the risk, and
averaging over several splits yields a {\em cross-validation}
estimate of the risk.
As will be shown in Sections~\ref{sec.def.simple2cross}
and~\ref{subsec.classical.examples}, various splitting strategies
lead to various CV estimates of the risk.

%%%avantages generaux, universalite
%
The major interest of CV lies in the universality of the data
splitting heuristics, which only assumes that data are identically
distributed and the training and validation samples are independent,
two assumptions which can even be relaxed (see
Section~\ref{sec.dependent}).
Therefore, CV can be applied to (almost) any algorithm in (almost)
any framework, for instance regression \citep{Sto:1974,Gei:1975},
density estimation \citep{Rude82,Ston84} and classification
\citep{DeWa79,Bar_Bou_Lug:2002}, among many others.
On the contrary, most other model selection procedures (see
Section~\ref{sec.modselproc}) are specific to a framework: For
instance, $C_p$ \citep{Mal:1973} is specific to least-squares
regression.

\subsection{From validation to cross-validation} \label{sec.def.simple2cross}
In this section, the hold-out (or validation) estimator of the risk
is defined, leading to a general definition of CV.

\subsubsection{Hold-out} \label{sec.def.simple2cross.ho}
The {\em hold-out} \citep{DeWa79} or (simple) {\em validation} relies
on a single split of data.
Formally, let $\It$ be a non-empty proper subset of $\set{1 , \ldots,
n}$, that is, such that both $\It$ and its complement $\Iv =
\paren{\It}^c = \set{1 , \ldots, n} \backslash \It$ are non-empty.
The {\em hold-out} estimator of the risk of $\A(D_n)$ with training
set $\It$ is defined by
\begin{equation} \label{def.validation}
\Loshval \paren{ \A ; D_n ; \It } \egaldef \frac{1}{n_v} \sum_{i \in
\Dv} \gamma\paren{ \A(\Dt) ; (X_i,Y_i)} \enspace ,
\end{equation}
where $\Dt \egaldef (\xi_i )_{i \in \It}$ is the {\em training
sample}, of size $n_t = \card(\It)$, and $\Dv \egaldef (\xi_i )_{i
\in \Iv}$ is the {\em validation sample}, of size $n_v = n-n_t$;
$\Iv$ is called the validation set.
The question of choosing $n_t$, and $\It$ given its cardinality
$n_t$, is discussed in the rest of this survey.

\subsubsection{General definition of cross-validation} \label{sec.def.simple2cross.CVgal}
A general description of the CV strategy has been given by
\cite{Gei:1975}: In brief, CV consists in averaging several hold-out
estimators of the risk corresponding to different splits of the data.
Formally, let $B \geq 1$ be an integer and $\It_1, \ldots, \It_B$ be
a sequence of non-empty proper subsets of $\set{1 , \ldots, n}$.
The CV estimator of the risk of $\A(D_n)$ with training sets
$\paren{\It_j}_{1 \leq j \leq B}$ is defined by
\begin{equation} \label{def.vc.gal}
\Loshvc \paren{ \A ; D_n ; \paren{\It_j}_{1 \leq j \leq B} } \egaldef
\frac{1}{B} \sum_{j=1}^B \Loshval \paren{ \A ; D_n ; \It_j } \enspace
.
\end{equation}
All existing CV estimators of the risk are of the form
\eqref{def.vc.gal}, each one being uniquely determined by the way the
sequence $\paren{\It_j}_{1 \leq j \leq B}$ is chosen, that is, the
choice of the splitting scheme.

%%% CV with voting

Note that when CV is used in model selection for identification, an
alternative definition of CV was proposed by
\cite{Yan:2006,Yan:2007b} and called {\em CV with voting} (CV-v).
When two algorithms $\A_1$ and $\A_2$ are compared, $\A_1$ is
selected by CV-v if and only if $ \Loshval ( \A_1  ; D_n ; \It_j ) <
\Loshval ( \A_2  ; D_n ; \It_j )$ for a majority of the splits $j =
1, \ldots, B$.
By contrast, CV procedures of the form \eqref{def.vc.gal} can be
called ``CV with averaging'' (CV-a), since the estimates of the risk
of the algorithms are averaged before their comparison.

\subsection{Classical examples}\label{subsec.classical.examples}
Most classical CV estimators split the data with a fixed size $n_t$
of the training set, that is, $\card(\It_j) \approx n_t$ for every
$j$. The question of choosing $n_t$ is discussed extensively in the
rest of this survey.
In this subsection, several CV estimators are defined. Two main
categories of splitting schemes can be distinguished, given $n_t$:
exhaustive data splitting, that is considering all training sets
$\It$ of size $n_t$, and partial data splitting.

\subsubsection{Exhaustive data splitting} \label{subsec.classical.examples.exhaust}
\paragraph{Leave-one-out} \citep[LOO,][]{Sto:1974,All:1974,Gei:1975} is the most classical exhaustive CV procedure, corresponding to the choice $n_t = n-1 \,$: Each data point is successively ``left out'' from the sample and used for validation.
Formally, LOO is defined by \eqref{def.vc.gal} with $B=n$ and $\It_j
= \set{j}^c$ for $j=1, \ldots, n \,$:
\begin{equation} \label{def.loo}
\Loshloo \paren{ \A ; D_n }
% = \frac{1}{n} \sum_{j=1}^n \Loshval \paren{ \A ; D_n ; \It_j }
= \frac{1}{n} \sum_{j=1}^n \gamma\paren{ \A \paren{D_n^{(-j)}} ;
\xi_j }
\end{equation}
where $D_n^{(-j)} = \paren{\xi_i}_{i \neq j}$.
%
% The LOO computation requires an algorithmic complexity $n$ times that
% of algorithm $\A$.
%
The name LOO can be traced back to papers by \cite{Pic_Coo:1984} and
by \cite{Bre_Spe:1992},
%\cite{Sha:1993} and by \cite{Zha:1993},
but LOO has several other names in the literature, such as {\em
delete-one CV} \citep[see][]{KCLi:1987}, {\em ordinary CV}
\citep{Sto:1974,Bur:1989}, or even only {\em CV}
\citep{Efr:1983,KCLi:1987}.

\paragraph{Leave-$p$-out} \citep[LPO,][]{Sha:1993} with $p \in \set{1, \ldots, n}$ is the exhaustive CV with $n_t = n-p\,$: every possible set of $p$ data points are successively ``left out'' from the sample and used for validation.
Therefore, LPO is defined by \eqref{def.vc.gal} with $B = {n\choose
p}$ and $\sparen{ \It_j}_{1 \leq j \leq B}$ are all the subsets of
$\set{1, \ldots, n}$ of size $p$.
LPO is also called {\em delete-$p$ CV} or {\em delete-$p$ multifold
CV} \citep{Zha:1993}.
Note that LPO with $p=1$ is LOO.

\subsubsection{Partial data splitting} \label{subsec.classical.examples.partial}
Considering ${n \choose p}$ training sets can be computationally
intractable, even for small $p$, so that partial data splitting
methods have been proposed.

\paragraph{{\it V}-fold CV} (VFCV) with $V \in \set{1, \ldots, n}$ was introduced by \cite{Gei:1975} as an
alternative to the computationally expensive LOO \citep[see also][for
instance]{Bre_etal:1984}.
VFCV relies on a preliminary partitioning of the data into $V$
subsamples of approximately equal cardinality $n/V$; each of these
subsamples successively plays the role of validation sample.
Formally, let $A_1, \ldots, A_V$ be some partition of $\set{1,
\ldots, n}$ with $\card\paren{A_j}\approx n/V$.
Then, the VFCV estimator of the risk of $\A$ is defined by
\eqref{def.vc.gal} with $B=V$ and $\It_j = A_j^c$ for $j = 1, \ldots,
B$, that is,
\begin{align} \label{def.vfcv}
\hspace*{-.15cm} \Loshvf \paren{ \ERM ; D_n ; \paren{ A_j }_{1 \leq j
\leq V}}
% &= \frac{1}{V} \sum_{j=1}^V \Loshval \paren{ \ERM ; D_n ; A_j^c } \\ \notag &
= \frac{1}{V} \sum_{j=1}^V \croch{ \frac{1}{\card(A_j)} \sum_{i \in
A_j} \gamma\paren{ \ERM \paren{D_n^{(-A_j)}} ; \xi_i } }
\end{align}
where $D_n^{(-A_j)} = \paren{\xi_i}_{i \in A_j^c}$.
By construction, the algorithmic complexity of VFCV is only $V$ times
that of training $\A$ with $n-n/V$ data points, which is much less
than LOO or LPO if $V \ll n$.
Note that VFCV with $V=n$ is LOO.

\paragraph{Balanced Incomplete CV} \citep[BICV,][]{Sha:1993} can be seen as an alternative to VFCV well-suited for small training sample sizes $n_t$.
Indeed, BICV is defined by \eqref{def.vc.gal} with training sets
$\paren{ A^c }_{A \in \mathcal{T}}$, where $\mathcal{T}$ is a
balanced incomplete block designs \citep[BIBD,][]{Joh:1971}, that is,
a collection of $B>0$ subsets of $\set{1, \ldots, n}$ of size
$n_v=n-n_t$ such that:
\begin{enumerate}
\item $\card\set{A \in \mathcal{T} \telque k \in \A}$ does not depend on $k \in \set{1, \ldots, n}$.
\item $\card\set{A \in \mathcal{T} \telque k,\ell \in \A}$ does not depend on $k \neq \ell \in \set{1, \ldots, n}$.
\end{enumerate}

The idea of BICV is to give to each data point (and each pair of data
points) the same role in the training and validation tasks.
Note that VFCV relies on a similar idea, since the set of training
sample indices used by VFCV satisfy the first property and almost the
second one: Pairs $(k,\ell)$ belonging to the same $A_j$ appear in
one validation set more than other pairs.

\paragraph{Repeated learning-testing} (RLT) was introduced by \cite{Bre_etal:1984} and further studied by
\cite{Bur:1989} and by \cite{Zha:1993} for instance.
The RLT estimator of the risk of $\A$ is defined by
\eqref{def.vc.gal} with any $B>0$ and $\sparen{\It_j}_{1 \leq j \leq
B}$ are $B$ different subsets of $\set{1, \ldots, n}$, chosen
randomly and independently from the data.
RLT can be seen as an approximation to LPO with $p = n - n_t$, with
which it coincides when $B = {n \choose p}$.

\paragraph{Monte-Carlo CV} \citep[MCCV,][]{Pic_Coo:1984} is very close to RLT: $B$ independent subsets of $\set{1, \ldots, n}$ are randomly drawn, with uniform distribution among subsets of size $n_t$.
The only difference with RLT is that MCCV allows the same split to be
chosen several times.

%%%

\subsubsection{Other cross-validation-like risk estimators} \label{sec.def.classex.other}
Several procedures have been introduced which are close to, or based
on CV. Most of them aim at fixing an observed drawback of CV.

\paragraph{Bias-corrected} versions of VFCV and RLT risk estimators have been proposed by \cite{Bur:1989,Bur:1990}, and a closely related penalization procedure called $V$-fold penalization has been defined by \cite{Arl:2008a}, see Section~\ref{sec.riskestim.bias.correction} for details.

\paragraph{Generalized CV} \citep[GCV,][]{Cra_Wah:1979} was introduced as a rotation-invariant version of LOO in least-squares regression, for estimating the risk of a linear estimator $\ERM = M {\bf Y}$ where ${\bf Y} = (Y_i)_{1 \leq i \leq n} \in \R^n$ and $M$ is an $n \times n$ matrix independent from ${\bf Y}$:
\begin{align*}
\crit_{\mathrm{GCV}} (M,{\bf Y}) \egaldef \frac{ n^{-1} \norm{{\bf Y}
- M {\bf Y}}^2 } { \paren{ 1 - n^{-1}  \tr(M)}^2} \quad \mbox{where}
\quad \forall t \in \R^n, \, \norm{t}^2 = \sum_{i=1}^n t_i^2 \enspace
.
\end{align*}
GCV is actually closer to $C_L$ \citep{Mal:1973} than to CV, since
GCV can be seen as an approximation to $C_L$ with a particular
estimator of the variance \citep{Efr:1986}. The efficiency of GCV has
been proved in various frameworks, in particular by
\cite{Li85,KCLi:1987} and by \cite{Cao_Gol:2006}.

\paragraph{Analytic Approximation}
When CV is used for selecting among linear models, \cite{Sha:1993}
proposed an analytic approximation to LPO with $p \sim n$, which is
called APCV.

\paragraph{LOO bootstrap and .632 bootstrap}
The bootstrap is often used for stabilizing an estimator or an
algorithm, replacing $\A(D_n)$ by the average of $\A(D_n^{\star})$
over several bootstrap resamples $D_n^{\star}$.
This idea was applied by \cite{Efr:1983} to the LOO estimator of the
risk, leading to the {\em LOO bootstrap}.
Noting that the LOO bootstrap was biased, \cite{Efr:1983} gave a
heuristic argument leading to the {\em $.632$ bootstrap} estimator of
the risk, later modified into the {\em $.632+$ bootstrap} by
\cite{Efr_Tib:1997}.
The main drawback of these procedures is the weakness of their
theoretical justifications. Only empirical studies have supported the
good behaviour of $.632+ $ bootstrap
\citep{Efr_Tib:1997,Mol_Sim_Pfe:2005}.

\subsection{Historical remarks}
Simple validation or hold-out was the first CV-like procedure. It was
introduced in the psychology area \citep{Lars31} from the need for a
reliable alternative to the {\em resubstitution error}, as
illustrated by \cite{AnAB72}. The hold-out was used by \cite{Herz69}
for assessing the quality of
predictors. % {\bf (???) by a theoretical and numerical study (???)}.
The problem of choosing the training set was first considered by
\cite{Sto:1974}, where ``controllable'' and ``uncontrollable'' data
splits were distinguished; an instance of uncontrollable division can
be found in the book by \cite{Simo71}.

A primitive LOO procedure was used by \cite{Hill66} and by
\cite{LaMi68} for evaluating the error rate of a prediction rule, and
a primitive formulation of LOO can be found in a paper by
\cite{MoTu68}. Nevertheless, LOO was actually introduced
independently by \cite{Sto:1974}, by \cite{All:1974} and by
\cite{Gei:1975}.
The relationship between LOO and the jackknife \citep{Que:1949},
which both rely on the idea of removing one observation from the
sample, has been discussed by \cite{Sto:1974} for instance.

The hold-out and CV were originally used only for estimating the risk
of an algorithm. The idea of using CV for model selection arose in
the discussion of a paper by \cite{EfMo73} and in a paper by
\cite{Geis74}. The first author to study LOO as a model selection
procedure was \cite{Sto:1974}, who proposed to use LOO again for
estimating the risk of the selected model.
%

%%%%%%%%%%%%%%%%%%%%%%%%%%%%%%%%%%%%%%%%%%%%%%%%%%%%%%%%%%%%%%%%%%%%
%%%%%%%%%%%%%%%%%%%%%%%%%%%%%%%%%%%%%%%%%%%%%%%%%%%%%%%%%%%%%%%%%%%%
%%%%%%%%%%%%%%%%%%%%%%%%%%%%%%%%%%%%%%%%%%%%%%%%%%%%%%%%%%%%%%%%%%%%
%%%%%%%%%%%%%%%%%%%%%%%%%%%%%%%%%%%%%%%%%%%%%%%%%%%%%%%%%%%%%%%%%%%%
%%%%%%%%%%%%%%%% Risk estimation
%\input{riskestim.tex}% Estimation of the risk for a single model

\section{Statistical properties of cross-validation
estimators of the risk} \label{sec.riskestim}
Understanding the behaviour of CV for model selection, which is the
purpose of this survey, requires first to analyze the performances of
CV as an estimator of the risk of a single algorithm.
Two main properties of CV estimators of the risk are of particular
interest: their bias, and their variance.

\subsection{Bias} \label{sec.riskestim.bias}

Dealing with the bias incurred by CV estimates can be made by two
strategies: evaluating the amount of bias in order to choose the
least biased CV procedure, or correcting for this bias.

\subsubsection{Theoretical assessment of the bias}
The independence of the training and the validation samples imply
that for every algorithm $\A$ and any $\It \subset \set{1, \ldots,
n}$ with cardinality $n_t$,
\[ \E\croch{ \Loshval \paren{ \A ; D_n ; \It } }
=  \E\croch{\gamma\paren{ \A \paren{\Dt} ; \xi }} = \E \croch{ \LosP{
\A\paren{D_{n_t}} }{P} } \enspace .
\]
Therefore, assuming that $\card(\It_j) = n_t$ for $j=1, \ldots, B$,
the expectation of the CV estimator of the risk only depends on $n_t
\,$:
\begin{equation} \label{eq.val.esperance}
\E\croch{ \Loshvc \paren{ \A ; D_n ; \paren{\It_j}_{1 \leq j \leq B}
} } = \E \croch{ \LosP{ \A\paren{D_{n_t}} }{P} } \enspace .
\end{equation}
In particular \eqref{eq.val.esperance} shows that the bias of the CV
estimator of the risk of $\A$ is the difference between the risks of
$\A$, computed respectively with $n_t$ and $n$ data points.
Since $n_t < n$, the bias of CV is usually nonnegative, which can be
proved rigorously when the risk of $\A$ is a decreasing function of
$n$, that is, when $\A$ is a smart rule; note however that a
classical algorithm such as 1-nearest-neighbour in classification is
not smart \citep[Section~6.8]{Dev_Gyo_Lug:1996}.
Similarly, the bias of CV tends to decrease with $n_t$, which is
rigorously true if $\A$ is smart.

More precisely, \eqref{eq.val.esperance} has led to several results
on the bias of CV, which can be split into three main categories:
asymptotic results ($\A$ is fixed and the sample size $n$ tends to
infinity), non-asymptotic results (where $\A$ is allowed to make use
of a number of parameters growing with $n$, say $n^{1/2}$, as often
in model selection), and empirical results.
They are listed below by statistical framework.

\paragraph{Regression}
The general behaviour of the bias of CV (positive, decreasing with
$n_t$) is confirmed by several papers and for several CV estimators.
For LPO, non-asymptotic expressions of its bias were proved by
\cite{Cel:2008:phd} for projection estimators, and by \cite{ArCe09}
for regressograms and kernels estimators when the design is fixed.
For VFCV and RLT, an asymptotic expansion of their bias was yielded
by \cite{Bur:1989} for least-squares estimators in linear regression,
and extended to spline smoothing \citep{Bur:1990}.
Note finally that \cite{Efr:1986} proved non-asymptotic analytic
expressions of the expectations of the LOO and GCV estimators of the
risk in regression with binary data \citep[see also][for some
explicit calculations]{Efr:1983}.

\paragraph{Density estimation}
shows a similar picture.
Non-asymptotic expressions for the bias of LPO estimators for kernel
and projection estimators with the quadratic risk were proved by
\cite{CeRo08} and by \cite{Celi08}.
Asymptotic expansions of the bias of the LOO estimator for histograms
and kernel estimators were previously proved by \cite{Rude82}; see
\cite{Bowm84} for simulations.
\cite{Hall87} derived similar results with the log-likelihood
contrast for kernel estimators, and related the performance of LOO to
the interaction between the kernel and the tails of the target
density $s$.

\paragraph{Classification}
For the simple problem of discriminating between two populations with
shifted distributions, \cite{Dav_Hal:1992} compared the asymptotical
bias of LOO and bootstrap, showing the superiority of the LOO when
the shift size is $n^{-1/2} \, $: As $n$ tends to infinity, the bias
of LOO stays of order $n^{-1}$, whereas that of bootstrap worsens to
the order $n^{-1/2}$.
On realistic synthetic and real biological data,
\cite{Mol_Sim_Pfe:2005} compared the bias of LOO, VFCV and .632+
bootstrap: The bias decreases with $n_t$, and is generally minimal
for LOO. Nevertheless, the $10$-fold CV bias is nearly minimal
uniformly over their experiments. In the same experiments, .632+
bootstrap exhibits the smallest bias for moderate sample sizes and
small signal-to-noise ratios, but a much larger bias otherwise.

\paragraph{CV-calibrated algorithms}
When a family of algorithm $\paren{\A_{\lambda}}_{\lambda \in
\Lambda}$ is given, and $\lh$ is chosen by minimizing
$\Loshvc(\A_{\lambda};D_n)$ over $\lambda$, $\Loshvc(\A_{\lh};D_n)$
is biased for estimating the risk of $\A_{\lh}(D_n)$, as reported
from simulation experiments by \cite{Sto:1974} for the LOO, and by
\cite{JoKM00} for VFCV in the variable selection setting.
This bias is of different nature compared to the previous frameworks.
Indeed, $\Loshvc(\A_{\lh},D_n)$ is biased simply because $\lh$ was
chosen using the same data as $\Loshvc(\A_{\lambda},D_n)$. This
phenomenon is similar to the optimism of $\LosP{\sh(D_n)}{P_n}$ as an
estimator of the loss of $\sh(D_n)$.
The correct way of estimating the risk of $\A_{\lh}(D_n)$ with CV is
to consider the full algorithm $\A': D_n \flapp \A_{\lh(D_n)}(D_n)$,
and then to compute $\Loshvc\paren{\A';D_n}$.
The resulting procedure is called ``double cross'' by
\cite{Sto:1974}.

\subsubsection{Correction of the bias}
\label{sec.riskestim.bias.correction}
An alternative to choosing the CV estimator with the smallest bias is
to correct for the bias of the CV estimator of the risk.
\cite{Bur:1989,Bur:1990} proposed a corrected VFCV estimator, defined
by
\begin{align*}
\Loshnom{corrVF} (\A; D_n) &= \Loshvf \paren{ \ERM ; D_n} +
\LosP{\A(D_n)}{P_n} - \frac{1}{V} \sum_{j=1}^V
\LosP{\A(D_n^{(-A_j)})}{P_n} \enspace ,
%\\ \Loshnom{corrRLT}  (\A; D_n) &= \Loshnom{RLT}  (\A; D_n) + \LosP{\A(D_n)}{P_n} - \frac{1}{B} \sum_{j=1}^B \LosP{\A(D_{n,j}^{(t)})}{P_n}
\end{align*}
and a corrected RLT estimator was defined similarly.
Both estimators have been proved to be asymptotically unbiased for
least-squares estimators in linear regression.

When the $A_j$s have exactly the same size $n/V$, the corrected VFCV
criterion is equal to the sum of the empirical risk and the $V$-fold
penalty \citep{Arl:2008a}, defined by
\[
\penVF (\A;D_n) = \frac{V-1}{V} \sum_{j=1}^V \croch{
\LosP{\A(D_n^{(-A_j)})}{P_n} - \LosP{\A(D_n^{(-A_j)})}{P_n^{(-A_j)}}
} \enspace .
\]
The $V$-fold penalized criterion was proved to be (almost) unbiased
in the non-asymptotic framework for regressogram estimators.

\subsection{Variance} \label{sec.riskestim.var}
CV estimators of the risk using training sets of the same size $n_t$
have the same bias, but they still behave quite differently; their
variance $\var(\Loshvc(\A;D_n;(\It_j)_{1 \leq j \leq B}))$ captures
most of the information to explain these differences.
\subsubsection{Variability factors}
Assume that $\card(\It_j) = n_t$ for every $j$. The variance of CV
results from the combination of several factors, in particular
$(n_t,n_v)$ and $B$.
%
%

%%%% influence de n_t et $n_v$
%
\paragraph{Influence of $(n_t,n_v)$}
Let us consider the hold-out estimator of the risk. Following in
particular \cite{NaBe03},
\begin{align}
& \qquad \notag \var\croch{ \Loshval\paren{\A;D_n;\It}} \\ \notag
&= \E\croch{\var\paren{ \LosP{\A(\Dt)}{P_n^{(v)}} \sachant \Dt }} + \var\croch{ \LosP{\A(D_{n_t})}{P} } \\
&= \frac{1}{n_v}\E\croch{\var\paren{\gamma\paren{\sh,\xi}\sachant \sh
= \A(\Dt)}} + \var\croch{ \LosP{\A(D_{n_t})}{P} }
\label{exp.variance}  \enspace .
\end{align}
The first term, proportional to $1/n_v$, shows that more data for
validation decreases the variance of $\Loshval$, because it yields a
better estimator of $\LosP{\A(\Dt)}{P}$.
The second term shows that the variance of $\Loshval$ also depends on
the distribution of $\LosP{\A(\Dt)}{P}$ around its expectation; in
particular, it strongly depends on the {\em stability} of $\A$.

\paragraph{Stability and variance}
When $\A$ is unstable, $\Loshloo\paren{\A}$ has often been pointed
out as a variable estimator \citep[Section~7.10,
][]{Has_Tib_Fri:2001,Bre:1996}.
Conversely, this trend disappears when $\A$ is stable, as noticed by
\cite{Mol_Sim_Pfe:2005} from a simulation experiment.

The relation between the stability of $\A$ and the variance of
$\Loshvc\paren{\A}$ was pointed out by \cite{DeWa79} in
classification, through upper bounds on the variance of
$\Loshloo\paren{\A}$.
\cite{BoEl02} extended these results to the regression setting, and
proved upper bounds on the maximal upward deviation of
$\Loshloo\paren{\A}$.

Note finally that several approaches based on the bootstrap have been
proposed for reducing the variance of $\Loshloo\paren{\A}$, such as
LOO bootstrap, .632 bootstrap and .632+ bootstrap \citep{Efr:1983};
see also Section~\ref{sec.def.classex.other}.

%%%% influence de B
%
\paragraph{Partial splitting and variance}
When $(n_t,n_v)$ is fixed, the variability of CV tends to be larger
for partial data splitting methods than for LPO. Indeed, having to
choose $B < {n \choose n_t}$ subsets $(\It_j)_{1 \leq j \leq B}$ of
$\set{1, \ldots, n}$, usually randomly, induces an additional
variability compared to $\Loshlpo$ with $p = n - n_t$.
In the case of MCCV, this variability decreases like $B^{-1}$ since
the $\It_j$ are chosen independently. The dependence on $B$ is
slightly different for other CV estimators such as RLT or VFCV,
because the $\It_j$ are not independent.
In particular, it is maximal for the hold-out, and minimal (null) for
LOO (if $n_t = n-1$) and LPO (with $p = n-n_t$).

Note that the dependence on $V$ for VFCV is more complex to evaluate,
since $B$, $n_t$, and $n_v$ simultaneously vary with $V$.
Nevertheless, a non-asymptotic theoretical quantification of this
additional variability of VFCV has been obtained by \cite{CeRo08} in
the density estimation framework \citep[see also empirical
considerations by][]{JoKM00}.

\subsubsection{Theoretical assessment of the variance}
Understanding precisely how $\var(\Loshvc(\A))$ depends on the
splitting scheme is complex in general, since $n_t$ and $n_v$ have a
fixed sum $n$, and the number of splits $B$ is generally linked with
$n_t$ (for instance, for LPO and VFCV).
Furthermore, the variance of CV behaves quite differently in
different frameworks, depending in particular on the stability of
$\A$. The consequence is that contradictory results have been
obtained in different frameworks, in particular on the value of $V$
for which the VFCV estimator of the risk has a minimal variance
\citep[Section~7.10]{Bur:1989,Has_Tib_Fri:2001}.
Despite the difficulty of the problem, the variance of several CV
estimators of the risk has been assessed in several frameworks, as
detailed below.

\paragraph{Regression}
In the linear regression setting, \cite{Bur:1989} yielded asymptotic
expansions of the variance of the VFCV and RLT estimators of the risk
with homoscedastic data.
The variance of RLT decreases with $B$, and in the case of VFCV, in a
particular setting,
\begin{align*}
\var\paren{\Loshvf(\A)} = \frac{2 \sigma^2}{n} + \frac{4
\sigma^4}{n^{2}} \croch{ 4 +
\frac{4}{V-1}+\frac{2}{(V-1)^2}+\frac{1}{(V-1)^3}}+o\paren{n^{-2}}
\enspace .
\end{align*}
The asymptotical variance of the VFCV estimator of the risk decreases
with $V$, implying that LOO asymptotically has the minimal variance.

Non-asymptotic closed-form formulas of the variance of the LPO
estimator of the risk have been proved by \cite{Cel:2008:phd} in
regression, for projection and kernel estimators for instance.
On the variance of RLT in the regression setting, see the asymptotic
results of \cite{Gir:1998} for Nadaraya-Watson kernel estimators, as
well as the non-asymptotic computations and simulation experiments by
\cite{NaBe03} with several learning algorithms.

\paragraph{Density estimation}
Non-asymptotic closed-form formulas of the variance of the LPO
estimator of the risk have been proved by \cite{CeRo08} and by
\cite{Celi08} for projection and kernel estimators.
In particular, the dependence of the variance of $\Loshnom{LPO}$ on
$p$ has been quantified explicitly for histogram and kernel
estimators by \cite{CeRo08}.

\paragraph{Classification}
For the simple problem of discriminating between two populations with
shifted distributions, \cite{Dav_Hal:1992} showed that the gap
between asymptotic variances of LOO and bootstrap becomes larger when
data are noisier.
\cite{NaBe03} made non-asymptotic computations and simulation
experiments with several learning algorithms.
\cite{Has_Tib_Fri:2001} empirically showed that VFCV has a minimal
variance for some $2<V<n$, whereas LOO usually has a large variance;
this fact certainly depends on the stability of the algorithm
considered, as showed by simulation experiments by
\cite{Mol_Sim_Pfe:2005}.

\subsubsection{Estimation of the variance} \label{sec.riskestim.var.estim}
There is no universal---valid under all distributions---unbiased
estimator of the variance of RLT \citep{NaBe03} and VFCV estimators
\citep{Ben_Gra:2004}.
In particular, \cite{Ben_Gra:2004} recommend the use of variance
estimators taking into account the correlation structure between test
errors; otherwise, the variance of CV can be strongly underestimated.

Despite these negative results, (biased) estimators of the variance
of $\Loshvc$ have been proposed by \cite{NaBe03}, by
\cite{Ben_Gra:2004} and by \cite{Mar_etal:2005}, and tested in
simulation experiments in regression and classification.
Furthermore, in the framework of density estimation with histograms,
\cite{CeRo08} proposed an estimator of the variance of the LPO risk
estimator. Its accuracy is assessed by a concentration inequality.
These results have recently been extended to projection estimators by
\cite{Celi08}.

%%%%%%%%%%%%%%%%%%%%%%%%%%%%%%%%%%%%%%%%%%%%%%%%%%%%%%%%%%%%%%%%%%%%%
%%%%%%%%%%%%%%%%%%%%%%%%%%%%%%%%%%%%%%%%%%%%%%%%%%%%%%%%%%%%%%%%%%%%%
%%%%%%%%%%%%%%%%%%%%%%%%%%%%%%%%%%%%%%%%%%%%%%%%%%%%%%%%%%%%%%%%%%%%%
%%%%%%%%%%%%%%%%%%%%%%%%%%%%%%%%%%%%%%%%%%%%%%%%%%%%%%%%%%%%%%%%%%%%%
%%%%%%%%%%%%%%%%%%% Efficiency
%\input{cvefficient.tex}% CV for efficient model selection

\section{Cross-validation for efficient model selection} \label{sec.cvefficient}
This section tackles the properties of CV procedures for model
selection when the goal is estimation (see
Section~\ref{sec.modsel.estim}).
\subsection{Relationship between risk estimation and model selection} \label{sec.cvefficient.link-risk-mod-sel}
As shown in Section~\ref{sec.modselproc.unbiased}, minimizing an
unbiased estimator of the risk leads to an efficient model selection
procedure.
One could conclude here that the best CV procedure for estimation is
the one with the smallest bias and variance (at least
asymptotically), for instance, LOO in the least-squares regression
framework \citep{Bur:1989}.

Nevertheless, the best CV estimator of the risk is not necessarily
the best model selection procedure. For instance, \cite{Bre_Spe:1992}
observed that uniformly over the models, the best risk estimator is
LOO, whereas 10-fold CV is more accurate for model selection.
Three main reasons for such a difference can be invoked.
First, the asymptotic framework ($\A$ fixed, $n \rightarrow \infty$)
may not apply to models close to the oracle, which typically has a
dimension growing with $n$ when $s$ does not belong to any model.
Second, as explained in Section~\ref{sec.modselproc.biased},
estimating the risk of each model with some bias can be beneficial
and compensate the effect of a large variance, in particular when the
signal-to-noise ratio is small.
Third, for model selection, what matters is not that every estimate
of the risk has small bias and variance, but more that
\[ \sign\paren{\crit(m_1) - \crit(m_2)} = \sign\paren{\LosP{\ERM_{m_1}}{P} - \LosP{\ERM_{m_2}}{P}} \]
with the largest probability for models $m_1,m_2$ near the oracle.

Therefore, specific studies are required to evaluate the performances
of the various CV procedures in terms of model selection efficiency.
In most frameworks, the model selection performance directly follows
from the properties of CV as an estimator of the risk, but not
always.

\subsection{The global picture}\label{subsubsec.asymptotic.equiv.}
Let us start with the classification of model selection procedures
made by \cite{Sha:1997} in the linear regression framework, since it
gives a good idea of the performance of CV procedures for model
selection in general.
Typically, the efficiency of CV only depends on the asymptotics of
$n_t/n \,$:
\begin{itemize}
\item When $n_t \sim n$, CV is asymptotically equivalent to Mallows' $C_p$, hence asymptotically optimal.
\item When $n_t \sim \lambda n$ with $\lambda \in (0,1)$, CV is asymptotically equivalent to GIC$_{\kappa}$ with $\kappa = 1+\lambda^{-1}$, which is defined as AIC with a penalty multiplied by $\kappa/2$. Hence, such CV procedures are overpenalizing by a factor $(1+\lambda)/(2 \lambda) > 1$.
\end{itemize}
The above results have been proved by \cite{Sha:1997} for LPO
\citep[see also][for the LOO]{KCLi:1987}; they also hold for RLT when
$B \gg n^2$ since RLT is then equivalent to LPO \citep{Zha:1993}.

In a general statistical framework, the model selection performance
of MCCV, VFCV, LOO, LOO Bootstrap, and .632 bootstrap for selection
among minimum contrast estimators was studied in a series of papers
\citep{vdL_Dud:2003,vdL_Dud_Kel:2004,vdL_Dud_vdV:2006,vdV_Dud_vdL:2006};
these results apply in particular to least-squares regression and
density estimation.
It turns out that under mild conditions, an oracle-type inequality is
proved, showing that up to a multiplying factor $C_n \rightarrow 1$,
the risk of CV is smaller than the minimum of the risks of the models
with a sample size $n_t$. In particular, in most frameworks, this
implies the asymptotic optimality of CV as soon as $n_t \sim n$. When
$n_t \sim \lambda n$ with $\lambda \in (0,1)$, this naturally
generalizes Shao's results.

\subsection{Results in various frameworks}
This section gathers results about model selection performances of CV
when the goal is estimation, in various frameworks. Note that model
selection is considered here with a general meaning, including in
particular bandwidth choice for kernel estimators.

\paragraph{Regression}
First, the results of Section~\ref{subsubsec.asymptotic.equiv.}
suggest that CV is suboptimal when $n_t$ is not asymptotically
equivalent to $n$. This fact has been proved rigorously for VFCV when
$V = \grandO(1)$ with regressograms \citep{Arl:2008a}:
with large probability, the risk of the model selected by VFCV is
larger than $1 + \kappa(V)$ times the risk of the oracle, with
$\kappa(V)>0$ for every fixed $V$.
Note however that the best $V$ for VFCV is not the largest one in
every regression framework, as shown empirically in linear regression
\citep{Bre_Spe:1992,HeTs86}; \cite{Bre:1996} proposed to explain this
phenomenon by relating the stability of the candidate algorithms and
the model selection performance of LOO in various regression
frameworks.

%
%

%%% heteroscedastic, regressograms
%
Second, the ``universality'' of CV has been confirmed by showing that
it naturally adapts to heteroscedasticity of data when selecting
among regressograms.
Despite its suboptimality, VFCV with $V = \grandO(1)$ satisfies a
non-asymptotic oracle inequality with constant $C>1$
\citep{Arl:2008a}.
Furthermore, $V$-fold penalization (which often coincides with
corrected VFCV, see Section~\ref{sec.riskestim.bias.correction})
satisfies a non-asymptotic oracle inequality  with $C_n \rightarrow
1$ as $n \rightarrow +\infty$, both when $V = \grandO(1)$
\citep{Arl:2008a} and when $V=n$ \citep{Arlo08c}.
Note that $n$-fold penalization is very close to LOO, suggesting that
it is also asymptotically optimal with heteroscedastic data.
Simulation experiments in the context of change-point detection
confirmed that CV adapts well to heteroscedasticity, contrary to
usual model selection procedures in the same framework
\citep{ArCe09}.

%%%% Autres cadres
%
The performances of CV have also been assessed for other kinds of
estimators in regression.
For choosing the number of knots in spline smoothing, \cite{Bur:1990}
proved that corrected versions of VFCV and RLT are asymptotically
optimal provided $n/(B n_v) = \grandO(1)$.
Furthermore, in kernel regression, several CV methods have been
compared to GCV in kernel regression by \cite{Har_Hal_Mar:1988} and
by \cite{Gir:1998}; the conclusion is that GCV and related criteria
are computationally more efficient than MCCV or RLT, for a similar
statistical performance.

%
%

%%%%%%% divers resultats
%
Finally, note that asymptotic results about CV in regression have
been proved by \cite{Gyo_etal:2002}, and an oracle inequality with
constant $C>1$ has been proved by \cite{Wegk03} for the hold-out,
with least-squares estimators.

\paragraph{Density estimation}
CV performs similarly than in regression for selecting among
least-squares estimators \citep{vdL_Dud_Kel:2004}: It yields a risk
smaller than the minimum of the risk with a sample size $n_t$. In
particular, non-asymptotic oracle inequalities with constant $C>1$
have been proved by \cite{Cel:2008:phd} for the LPO when $p/n\in
[a,b]$, for some $0<a<b<1$.

The performance of CV for selecting the bandwidth of kernel density
estimators has been studied in several papers.
With the least-squares contrast, the efficiency of LOO was proved by
\cite{Hal:1983} and generalized to the multivariate framework by
\cite{Ston84}; an oracle inequality asymptotically leading to
efficiency was recently proved by \cite{Dale05}.
%
%
%
%
%%% kernel Kullback
%
%
With the Kullback-Leibler divergence, CV can suffer from troubles in
performing model selection \citep[see also][]{ScGr81,ChGW87}.
The influence of the tails of the target $s$ was studied by
\cite{Hall87}, who gave conditions under which CV is efficient and
the chosen bandwidth is optimal at first-order.

\paragraph{Classification}
In the framework of binary classification by intervals (that is, with
$\X=[0,1]$ and piecewise constant classifiers), \cite{KMNR97} proved
an oracle inequality for the hold-out.
Furthermore, empirical experiments show that CV yields (almost)
always the best performance, compared to deterministic penalties
\citep{KMNR97}.
On the contrary, simulation experiments by \cite{Bar_Bou_Lug:2002} in
the same setting showed that random penalties such as Rademacher
complexity and maximal discrepancy usually perform much better than
hold-out, which is shown to be more variable.

Nevertheless, the hold-out still enjoys quite good theoretical
properties: It was proved to adapt to the margin condition by
\cite{Bla_Mas:2006}, a property nearly unachievable with usual model
selection procedures \citep[see
also][Section~8.5]{Mas:2003:St-Flour}.
This suggests that CV procedures are naturally adaptive to several
unknown properties of  data in the statistical learning framework.

The performance of the LOO in binary classification was related to
the stability of the candidate algorithms by \cite{KeRo99}; they
proved oracle-type inequalities called ``sanity-check bounds'',
describing the worst-case performance of LOO \citep[see
also][]{BoEl02}.

An experimental comparison of several CV methods and bootstrap-based
CV (in particular .632+ bootstrap) in classification can also be
found in papers by \cite{Efr:1986} and \cite{Efr_Tib:1997}.

%%%%%%%%%%%%%%%%%%%%%%%%%%%%%%%

%%%%%%%%%%%%%%%%%%%%%%%%%%%%%%%%%%%%%%%%%%%%%%%%%%%%%%%%%%%%%%%%%%%%%%
%%%%%%%%%%%%%%%%%%%%%%%%%%%%%%%%%%%%%%%%%%%%%%%%%%%%%%%%%%%%%%%%%%%%%%
%%%%%%%%%%%%%%%%%%%%%%%%%%%%%%%%%%%%%%%%%%%%%%%%%%%%%%%%%%%%%%%%%%%%%%
%%%%%%%%%%%%%%%%%%%%%%%%%%%%%%%%%%%%%%%%%%%%%%%%%%%%%%%%%%%%%%%%%%%%%%
%%%%%%%%%%%%%% Consistency
%\input{cvconsistent.tex}% CV for consistent model selection

\section{Cross-validation for identification} \label{sec.cvconsistent}
Let us now focus on model selection when the goal is to identify the
``true model'' $S_{m_0}$, as described in
Section~\ref{sec.modsel.identif}.
In this framework, asymptotic optimality is replaced by (model)
consistency, that is,
\[ \Proba \paren{ \mh(D_n) = m_0 } \xrightarrow[n \rightarrow \infty]{} 1 \enspace .\]
Classical model selection procedures built for identification, such
as BIC, are described in Section~\ref{sec.modselproc.identif}.

\subsection{General conditions towards model consistency}

At first sight, it may seem strange to use CV for identification:
LOO, which is the pioneering CV procedure, is actually closely
related to the unbiased risk estimation principle, which is only
efficient when the goal is estimation. Furthermore, estimation and
identification are somehow contradictory goals, as explained in
Section~\ref{sec.modsel.estim-vs-identif}.

This intuition about inconsistency of some CV procedures is confirmed
by several theoretical results.
\cite{Sha:1993} proved that several CV methods are inconsistent for
variable selection in linear regression: LOO, LPO, and BICV when
$\liminf_{n \rightarrow \infty} (n_t/n) > 0$.
Even if these CV methods asymptotically select all the true variables
with probability~1, the probability that they select too much
variables does not tend to zero.
More generally, \cite{Sha:1997} proved that CV procedures behave
asymptotically like GIC$_{\lambda_n}$ with $\lambda_n = 1 + n/n_t$,
which leads to inconsistency as soon as $n / n_t = \grandO(1)$.

In the context of ordered variable selection in linear regression,
\cite{Zha:1993} computed the asymptotic value of the probability of
selecting the true model for several CV procedures. He also
numerically compared the values of this probability for the same CV
procedures in a specific example. For LPO with $p/n \rightarrow
\lambda \in ( 0 , 1)$ as $n$ tends to $+\infty$, $\Prob\paren{\mh =
m_0}$ increases with $\lambda$.
The result is slightly different for VFCV: $\Prob\paren{\mh = m_0}$
 increases with $V$ (hence, it is maximal for the LOO, which is the worst case of LPO).
The variability induced by the number $V$ of splits seems to be more
 important here than the bias of VFCV. Nevertheless, $\Prob\paren{\mh = m_0}$ is almost
 constant between $V=10$ and $V=n$, so that taking $V>10$ is not advised for computational reasons.

\medskip

These results suggest that if the training sample size $n_t$ is
negligible in front of $n$, then model consistency could be obtained.
This has been confirmed theoretically by \cite{Sha:1993,Sha:1997} for
the variable selection problem in linear regression: CV is consistent
when $n \gg n_t \rightarrow \infty$, in particular RLT, BICV (defined
in Section~\ref{subsec.classical.examples.partial}) and LPO with
$p=p_n \sim n$ and $n-p_n \rightarrow \infty$.

Therefore, when the goal is to identify the true model, a larger
proportion of the data should be put in the validation set in order
to improve the performance. This phenomenon is somewhat related to
the {\em cross-validation paradox} \citep{Yan:2006}.

\subsection{Refined analysis for the algorithm selection problem}

The behaviour of CV for identification is better understood by
considering a more general framework, where the goal is to select
among statistical algorithms the one with the fastest convergence
rate.
\cite{Yan:2006,Yan:2007b} considered this problem for two candidate
algorithms (or more generally any finite number of algorithms).
Let us mention here that \cite{Sto:1977b} considered a few specific
examples of this problem, and showed that LOO can be inconsistent for
choosing the best among two ``good'' estimators.

The conclusion of Yang's papers is that the sufficient condition on
$n_t$ for the consistency in selection of CV strongly depends on the
convergence rates $\paren{r_{n,i}}_{i=1,2}$ of the candidate
algorithms. Let us assume that $r_{n,1}$ and $r_{n,2}$ differ at
least by a multiplicative constant $C>1$.
Then, in the regression framework, if the risk of $\ERM_i$ is
measured by $\E \norm{\ERM_i - s}_2$, \cite{Yan:2007b} proved that
the hold-out, VFCV, RLT and LPO with voting (CV-v, see
Section~\ref{sec.def.simple2cross.CVgal}) are consistent in selection
if
\begin{equation} \label{eq.cond-consist.reg}
n_v, n_t \rightarrow \infty \quad \mbox{and} \quad \sqrt{n_v} \max_i
r_{n_t,i} \rightarrow \infty
%\quad \mbox{where} \quad n_v = n - n_t
\enspace ,
\end{equation}
under some conditions on $\norm{\ERM_i - s}_p$ for $p=2,4,\infty$.
In the classification framework, if the risk of $\ERM_i$ is measured
by $\Prob\paren{\ERM_i \neq s}$, \cite{Yan:2006} proved the same
consistency result for CV-v under the condition
\begin{equation} \label{eq.cond-consist.classif}
n_v, n_t \rightarrow \infty \quad \mbox{and} \quad \frac{ n_v \max_i
r_{n_t,i}^2 } { s_{n_t}} \rightarrow \infty
%\quad \mbox{where} \quad n_v = n - n_t
\enspace ,
\end{equation}
where $s_n$ is the convergence rate of $\Prob\paren{\ERM_1(D_n) \neq
\ERM_2(D_n)}$.

Intuitively, consistency holds as soon as the uncertainty of each
estimate of the risk (roughly proportional to $n_v^{-1/2}$) is
negligible in front of the risk gap $\absj{r_{n_t,1} - r_{n_t,2}}$
(which is of the same order as $\max_i r_{n_t,i}$). This condition
holds either when at least one of the algorithms converges at a
non-parametric rate, or when $n_t \ll n$, which artificially widens
the risk gap.
Empirical results in the same direction were proved by
\cite{Die:1998} and by \cite{Alp:1999}, leading to the advice that
$V=2$ is the best choice when VFCV is used for comparing two learning
procedures. See also the results by \cite{NaBe03} about CV considered
as a testing procedure comparing two candidate algorithms.

The sufficient conditions \eqref{eq.cond-consist.reg} and
\eqref{eq.cond-consist.classif} can be simplified depending on
$\max_i r_{n,i}$, so that the ability of CV to distinguish between
two algorithms depends on their convergence rates.
On the one hand, if $\max_i r_{n,i} \propto n^{-1/2}$, then
\eqref{eq.cond-consist.reg} or \eqref{eq.cond-consist.classif} only
hold when $n_v \gg n_t \rightarrow \infty$ (under some conditions on
$s_n$ in classification). Therefore, the cross-validation paradox
holds for comparing algorithms converging at the parametric rate
(model selection when a true model exists being only a particular
case).
Note that possibly stronger conditions can be required in
classification where algorithms can converge at fast rates, between
$n^{-1}$ and $n^{-1/2}$.

On the other hand, \eqref{eq.cond-consist.reg} and
\eqref{eq.cond-consist.classif} are milder conditions when $\max_i
r_{n,i} \gg n^{-1/2}$: They are implied by $n_t / n_v = \grandO(1)$,
and they even allow $n_t \sim n$  (under some conditions on $s_n$ in
classification).
Therefore, non-parametric algorithms can be compared by more usual CV
procedures ($n_t>n/2$), even if LOO is still excluded by conditions
\eqref{eq.cond-consist.reg} and \eqref{eq.cond-consist.classif}.

\medskip

Note that according to a simulation experiments, CV with averaging
(that is, CV as usual) and CV with voting are equivalent at first but
not at second order, so that they can differ when $n$ is small
\citep{Yan:2007b}.

%

%%%%%%%%%%%%%%%%%%%%%%%%%%%%%%%%%%%%%%%%%%%%%%%%%%%%%%%%%%%%%%%%%%%
%%%%%%%%%%%%%%%%%%%%%%%%%%%%%%%%%%%%%%%%%%%%%%%%%%%%%%%%%%%%%%%%%%%
%%%%%%%%%%%%%%%%%%%%%%%%%%%%%%%%%%%%%%%%%%%%%%%%%%%%%%%%%%%%%%%%%%%
%%%%%%%%%%%%%%%%%%%%%%%%%%%%%%%%%%%%%%%%%%%%%%%%%%%%%%%%%%%%%%%%%%%
%%%%%%%%%%%%%%%% Specific frameworks
%\input{specific.tex}% Specificities of some frameworks

\section{Specificities of some frameworks} \label{sec.specific}

Originally, the CV principle has been proposed for {\em i.i.d.}
observations and usual contrasts such as least-squares and
log-likelihood. Therefore, CV procedures may have to be modified in
other specific frameworks, such as estimation in presence of outliers
or with dependent data.

\subsection{Density estimation}
In the density estimation framework, some specific modifications of
CV have been proposed.

First, \cite{Hal_Mar_Par:1992} defined the ``smoothed CV'', which
consists in pre-smoothing the data before using CV, an idea related
to the smoothed bootstrap.
They proved that smoothed CV yields an excellent asymptotical model
selection performance under various smoothness conditions on the
density.

Second, when the goal is to estimate the density at one point (and
not globally), \cite{Hal_Sch:1989} proposed a local version of CV and
proved its asymptotic optimality.

\subsection{Robustness to outliers}

In presence of outliers in regression, \cite{Leu:2005} studied how CV
must be modified to get both asymptotic efficiency and a consistent
bandwidth estimator \citep[see also][]{LeMW93}.

Two changes are possible to achieve robustness: Choosing a ``robust''
regressor, or choosing a robust loss-function.
In presence of outliers, classical CV with a non-robust loss function
has been shown to fail by \cite{Haerd84}.

\cite{Leu:2005} described a CV procedure based on robust losses like
$L^1$ and Huber's \citep{Hube64} ones.
The same strategy remains applicable to other setups like linear
models in \cite{RoFB97}.

\subsection{Time series and dependent observations} \label{sec.dependent}
As explained in Section~\ref{sec.def.philo}, CV is built upon the
heuristics that part of the sample (the validation set) can play the
role of {\em new data} with respect to the rest of the sample (the
training set).
``New'' means that the validation set is independent from the
training set with the same distribution.

Therefore, when data $\xi_1, \ldots, \xi_n$ are not independent, CV
must be modified, like other model selection procedures \citep[in
non-parametric regression with dependent data, see the review
by][]{Ops_Wan_Yan:2001}.

\medskip

Let us first consider the statistical framework of
Section~\ref{sec.intro} with $\xi_1, \ldots, \xi_n$ identically
distributed but not independent.
Then, when for instance data are positively correlated,
\cite{Har_Weh:1986} proved that CV overfits for choosing the
bandwidth of a kernel estimator in regression \citep[see
also][]{Chu_Mar:1991,Ops_Wan_Yan:2001}.

The main approach used in the literature for solving this issue is to
choose $\It$ and $\Iv$ such that $\min_{i \in \It, \, j \in \Iv}
\absj{i-j} > h > 0$, where $h$ controls the distance from which
observations $i$ and $j$ are independent.
For instance, the LOO can be changed into: $\Iv = \set{J}$ where $J$
is uniformly chosen in $\set{1, \ldots, n}$, and $\It = \set{1,
\ldots, J-h-1, J+h+1, \ldots, n}$, a method called ``modified CV'' by
\cite{Chu_Mar:1991} in the context of bandwidth selection.
%
%The main difference with usual CV procedures is that $\It \cup \Iv \neq \set{1, \ldots, n}$.
%
Then, for short range dependences, $\xi_i$ is almost independent from
$\xi_j$ when $\absj{i-j} > h$ is large enough, so that
$\paren{\xi_j}_{j \in \It}$ is almost independent from
$\paren{\xi_j}_{j \in \Iv}$.
Several asymptotic optimality results have been proved on modified
CV, for instance by \cite{Har_Vie:1990} for bandwidth choice in
kernel density estimation, when data are $\alpha$-mixing (hence, with
a short range dependence structure) and $h=h_n \rightarrow \infty$
``not too fast''.
Note that modified CV also enjoys some asymptotic optimality results
with long-range dependences, as proved by \cite{Hal_Lah_Pol:1995},
even if an alternative block bootstrap method seems more appropriate
in such a framework.

\medskip

Several alternatives to modified CV have also been proposed.
The ``$h$-block CV'' \citep{Bur_Cho_Nol:1994} is modified CV plus a
corrective term, similarly to the bias-corrected CV by
\cite{Bur:1989} (see Section~\ref{sec.riskestim.bias}).
Simulation experiments in several (short range) dependent frameworks
show that this corrective term matters when $h/n$ is not small, in
particular when $n$ is small.

The ``partitioned CV'' has been proposed by \cite{Chu_Mar:1991} for
bandwidth selection: An integer $g>0$ is chosen, a bandwidth $\lh_k$
is chosen by CV based upon the subsample $\paren{\xi_{k+g j}}_{j \geq
0}$ for each $k = 1, \ldots, g$, and the selected bandwidth is a
combination of $\sparen{\lh_k}$.

When a parametric model is available for the dependency structure,
\cite{Har:1994} proposed the ``time series CV''.

\medskip

An important framework where data often are dependent is time-series
analysis, in particular when the goal is to predict the next
observation $\xi_{n+1}$ from the past $\xi_1, \ldots, \xi_n$.
When data are stationary, $h$-block CV and similar approaches can be
used to deal with (short range) dependences. Nevertheless,
\cite{Bur_Nol:1992} proved in some specific framework that unaltered
CV is asymptotic optimal when $\xi_1, \ldots, \xi_n$ is a stationary
Markov process.
%
%%%

On the contrary, using CV for non-stationary time-series is a quite
difficult problem.
The only reasonable approach in general is the hold-out, that is,
$\It = \set{1, \ldots, m}$ and $\Iv = \set{m+1, \ldots, n}$ for some
deterministic $m$. Each model is first trained with $\paren{\xi_j}_{j
\in \It}$. Then, it is used for predicting successively $\xi_{m+1}$
from $\paren{\xi_j}_{j \leq m}$, $\xi_{m+2}$ from $\paren{\xi_j}_{j
\leq m+1}$, and so on. The model with the smallest average error for
predicting $\paren{\xi_j}_{j \in \Iv}$ from the past is chosen.

\subsection{Large number of models} \label{sec.specific.manymodels}
As mentioned in Section~\ref{sec.modselproc}, model selection
procedures estimating unbiasedly the risk of each model fail when, in
particular, the number of models grows exponentially with $n$
\citep{Bir_Mas:2006}.
Therefore, CV cannot be used directly, except maybe with $n_t \ll n$,
provided $n_t$ is well chosen \citep[see
Section~\ref{sec.cvefficient} and][Chapter~6]{Cel:2008:phd}.

For least-squares regression with homoscedastic data, \cite{Wegk03}
proposed to add to the hold-out estimator of the risk a penalty term
depending on the number of models. This method is proved to satisfy a
non-asymptotic oracle inequality with leading constant $C >1$.

Another general approach was proposed by \cite{ArCe09} in the context
of multiple change-point detection. The idea is to perform model
selection in two steps:
First, gather the models $\paren{S_m}_{\mM_n}$ into meta-models
$\sparen{\St_D}_{D \in \D_n}$, where $\D_n$ denotes a set of indices
such that $\card(\D_n)$ grows at most polynomially with $n$. Inside
each meta-model $\St_D = \bigcup_{\mM_n(D)} S_m$, $\ERM_D$ is chosen
from data by optimizing a given criterion, for instance the empirical
contrast $\LosP{t}{P_n}$, but other criteria can be used.
Second, CV is used for choosing among $\paren{\ERM_D}_{D \in \D_n}$.
Simulation experiments show this simple trick automatically takes
into account the cardinality of $\M_n$, even when data are
heteroscedastic, contrary to other model selection procedures built
for exponential collection of models which all assume
homoscedasticity of data.

%%%%%%%%%%%%%%%%%%%%%%%%%%%%%%%%%%%%%%%%%%%%%%%%%%%%%%%%%%%%%%%%%%%%%%%%
%%%%%%%%%%%%%%%%%%%%%%%%%%%%%%%%%%%%%%%%%%%%%%%%%%%%%%%%%%%%%%%%%%%%%%%%
%%%%%%%%%%%%%%%%%%%%%%%%%%%%%%%%%%%%%%%%%%%%%%%%%%%%%%%%%%%%%%%%%%%%%%%%
%%%%%%%%%%%%%%%%%%%%%%%%%%%%%%%%%%%%%%%%%%%%%%%%%%%%%%%%%%%%%%%%%%%%%%%%
%%%%%%%%%%%%%% Practical considerations
%\input{practice.tex}% Practical questions

\section{Closed-form formulas and fast computation} \label{sec.complex}

Resampling strategies, like CV, are known to be time consuming.
The naive implementation of CV has a computational complexity of $B$
times the complexity of training each algorithm $\A$, which is
usually  intractable for LPO, even with $p=1$. The computational cost
of VFCV or RLT can still be quite costly when $B>10$ in many
practical problems.
Nevertheless, closed-form formulas for CV estimators of the risk can
be obtained in several frameworks, which greatly decreases the
computational cost of CV.

\medskip

In density estimation, closed-form formulas have been originally
derived by \cite{Rude82} and by \cite{Bowm84} for the LOO risk
estimator of histograms and  kernel estimators. These results have
been recently extended by \cite{CeRo08} to the LPO risk estimator
with the quadratic loss. Similar results are more generally available
for projection estimators as settled by \cite{Celi08}.
Intuitively, such formulas can be obtained provided the number $N$ of
values taken by the $B={n\choose n_v}$ hold-out estimators of the
risk, corresponding to different data splittings, is at most
polynomial in the sample size.

\medskip

For least-squares estimators in linear regression, \cite{Zha:1993}
proved a closed-form formula for the LOO estimator of the risk.
Similar results have been obtained by \cite{Wahb75,Wahb77a}, and by
\cite{Cra_Wah:1979} in the spline smoothing context as well.
These papers led in particular to the definition of GCV (see
Section~\ref{sec.def.classex.other}) and related procedures, which
are often used instead of CV (with a naive implementation) because of
their small computational cost, as emphasized by \cite{Gir:1998}.

%%%%%

Closed-form formulas for the LPO estimator of the risk were also
obtained by \cite{Cel:2008:phd} in regression for kernel and
projection estimators, in particular for regressograms.
An important property of these closed-form formulas is their
additivity: For a regressogram associated to a partition
$(\Il)_{\lamm}$ of $\X$, the LPO estimator of the risk can be written
as a sum over $\lamm$ of terms which only depend on observations
$(X_i,Y_i)$ such that $X_i \in \Il$.
Therefore, dynamic programming \citep{BeDr62} can be used for
minimizing the LPO estimator of the risk over the set of partitions
of $\X$ in $D$ pieces. As an illustration, \cite{ArCe09} successfully
applied this strategy in the change-point detection framework.
Note that the same idea can be used with VFCV or RLT, but for a
larger computational cost since no closed-form formulas are available
for these CV methods.

\medskip

Finally, in frameworks where no closed-form formula can be proved,
some efficient algorithms exist for avoiding to recompute
$\Loshval\sparen{\A;D_n;\It_j}$ from scratch for each data splitting
$\It_j$.
These algorithms rely on updating formulas such as the ones by
\cite{Ripl96} for LOO in linear and quadratic discriminant analysis;
this approach makes LOO as expensive to compute as the empirical
risk.

Very similar formulas are also available for LOO and the $k$-nearest
neighbours estimator in classification \citep{DaMa08}.

%

%%%%%%%%%%%%%%%%%%%%%%%%%%%%%%%%%%%%%%%%%%%%%%%%%%%%%%%%%%%%%%%%%
%%%%%%%%%%%%%%%%%%%%%%%%%%%%%%%%%%%%%%%%%%%%%%%%%%%%%%%%%%%%%%%%%
%%%%%%%%%%%%%%%%%%%%%%%%%%%%%%%%%%%%%%%%%%%%%%%%%%%%%%%%%%%%%%%%%
%%%%%%%%%%%%%%%%%%%%%%%%%%%%%%%%%%%%%%%%%%%%%%%%%%%%%%%%%%%%%%%%%
%%%%%%%%%%%%%% Conclusion and guidelines
%\input{conclu.tex}% Conclusion

\section{Conclusion: which cross-validation method for which problem?} \label{sec.conclu}
This conclusion collects a few guidelines aiming at helping CV users,
first interpreting the results of CV, second appropriately using CV
in each specific problem.

\subsection{The general picture} \label{sec.conclu.gal}
Drawing a general conclusion on CV methods is an impossible task
because of the variety of frameworks where CV can be used, which
induces a variety of behaviors of CV. Nevertheless, we can still
point out the three main criteria to take into account for choosing a
CV method for a particular model selection problem:
\begin{itemize}
\item {\em Bias}: CV roughly estimates the risk of a model with a sample size $n_t < n$ (see Section~\ref{sec.riskestim.bias}). Usually, this implies that CV overestimates the variance term compared to the bias term in the bias-variance decomposition \eqref{eq.biais-var} with sample size $n$. \\
When the goal is estimation and the signal-to-noise ratio (SNR) is large, the smaller bias usually is the better, which is obtained by taking $n_t \sim n$. Otherwise, CV can be asymptotically suboptimal. Nevertheless, when the goal is estimation and the SNR is small, keeping a small upward bias for the variance term often improves the performance, which is obtained by taking $n_t \sim \kappa n$ with $\kappa \in (0,1)$. See Section~\ref{sec.cvefficient}. \\
When the goal is identification, a large bias is often needed, which
is obtained by taking $n_t \ll n$; depending on the framework, larger
values of $n_t$ can also lead to model consistency, see
Section~\ref{sec.cvconsistent}.
\item {\em Variability}: The variance of the CV estimator of the risk is usually a decreasing function of the number $B$ of splits, for a fixed training size.
When the number of splits is fixed, the variability of CV also
depends on the training sample size $n_t$. Usually, CV is more
variable when $n_t$ is closer to $n$.
However, when $B$ is linked with $n_t$ (as for VFCV or LPO), the
variability of CV must be quantified precisely, which has been done
in few frameworks. The only general conclusion on this point is that
the CV method with minimal variability seems strongly
framework-dependent, see Section~\ref{sec.riskestim.var} for details.
\item {\em Computational complexity}: Unless closed-form formulas or analytic approximations are available (see Section~\ref{sec.complex}), the complexity of CV is roughly proportional to the number of data splits: 1 for the hold-out, $V$ for VFCV, $B$ for RLT or MCCV, $n$ for LOO, and $\binom{n}{p}$ for LPO.
\end{itemize}
%
%\medskip
%
The optimal trade-off between these three factors can be different
for each problem, depending for instance on the computational
complexity of each estimator, on specificities of the framework
considered, and on the final user's trade-off between statistical
performance and computational cost.
Therefore, no ``optimal CV method'' can be pointed out before having
taken into account the final user's preferences.

Nevertheless, in density estimation, closed-form expressions of the
LPO estimator have been derived by \cite{CeRo08} with histograms and
kernel estimators, and by \cite{Celi08} for projection estimators.
These expressions allow to perform LPO without additional
computational cost, which reduces the aforementioned trade-off to the
easier bias-variability trade-off. In particular, \cite{CeRo08}
proposed to choose $p$ for LPO by minimizing a criterion defined as
the sum of a squared bias and a variance terms \citep[see
also][Chapter~9]{Pol_Rom_Wol:1999}.

\subsection{How the splits should be chosen?} \label{sec.choix-It}

For hold-out, VFCV, and RLT, an important question is to choose a
particular sequence of data splits.

First, should this step be random and independent from $D_n$, or take
into account some features of the problem or of the data?
It is often recommended to take into account the structure of data
when choosing the splits. If data are stratified, the proportions of
the different strata should (approximately) be the same in the sample
and in each training and validation sample.
Besides, the training samples should be chosen so that $\ERM_m(\Dt)$
is well defined for every training set; in the regressogram case,
this led \cite{Arl:2008a} and \cite{ArCe09} to choose carefully the
splitting scheme.
In supervised classification, practitioners usually choose the splits
so that the proportion of each class is the same in every validation
sample as in the sample.
Nevertheless, \cite{Bre_Spe:1992} made simulation experiments in
regression for comparing several splitting strategies. No significant
improvement was reported from taking into account the stratification
of data for choosing the splits.

Another question related to the choice of $(\It_j)_{1 \leq j \leq B}$
is whether the $\It_j$ should be independent (like MCCV), slighly
dependent (like RLT), or strongly dependent (like VFCV).
It seems intuitive that giving similar roles to all data points in
the $B$ ``training and validation tasks'' should yield more reliable
results as other methods. This intuition may explain why VFCV is much
more used than RLT or MCCV.
Similarly, \cite{Sha:1993} proposed a CV method called BICV, where
every point and pair of points appear in the same number of splits,
see Section~\ref{subsec.classical.examples.partial}.
Nevertheless, most recent theoretical results on the various CV
procedures are not accurate enough to distinguish which one may be
the best splitting strategy: This remains a widely open theoretical
question.

Note finally that the additional variability due to the choice of a
sequence of data splits was quantified empirically by \cite{JoKM00}
and theoretically by \cite{CeRo08} for VFCV.

\subsection{V-fold cross-validation} \label{sec.practice.VFCV}
VFCV is certainly the most popular CV procedure, in particular
because of its mild computational cost. Nevertheless, the question of
choosing $V$ remains widely open, even if indications can be given
towards an appropriate choice.

%%%%%% General
%
A specific feature of VFCV---as well as exhaustive strategies---is
that choosing $V$ uniquely determines the size of the training set
$n_t=n(V-1)/V$ and the number of splits $B=V$, hence the
computational cost.
Contradictory phenomena then occur.

On the one hand, the bias of VFCV decreases with $V$ since
$n_t=n(1-1/V)$ observations are used in the training set.
On the other hand, the variance of VFCV decreases with $V$ for small
values of $V$, whereas the LOO ($V=n$) is known to suffer from a high
variance in several frameworks such as classification or density
estimation. Note however that the variance of VFCV is minimal for
$V=n$ in some frameworks like linear regression (see
Section~\ref{sec.riskestim.var}).
Furthermore, estimating the variance of VFCV from data is a difficult
problem in general, see Section~\ref{sec.riskestim.var.estim}.

\medskip

When the goal of model selection is estimation, it is often reported
in the literature that the optimal $V$ is between $5$ and~$10$,
because the statistical performance does not increase much for larger
values of $V$, and averaging over 5 or~10 splits remains
computationally feasible \citep[][Section~7.10]{Has_Tib_Fri:2001}.
Even if this claim is clearly true for many problems, the conclusion
of this survey is that better statistical performance can sometimes
be obtained with other values of $V$, for instance depending on the
SNR value.

When the SNR is large, the asymptotic comparison of CV procedures
recalled in Section~\ref{subsubsec.asymptotic.equiv.} can be trusted:
LOO performs (nearly) unbiased risk estimation hence is
asymptotically optimal, whereas VFCV with $V=\grandO(1)$ is
suboptimal.
On the contrary, when the SNR is small, overpenalization can improve
the performance. Therefore, VFCV with $V<n$ can yield a smaller risk
than LOO thanks to its bias and despite its variance when $V$ is
small \citep[see simulation experiments by][]{Arl:2008a}.
Furthermore, other CV procedures like RLT can be interesting
alternatives to VFCV, since they allow to choose the bias (through
$n_t$) independently from $B$, which mainly governs the variance.
Another possible alternative is $V$-fold penalization, which is
related to corrected VFCV (see Section~\ref{sec.def.classex.other}).

%

%%%%%% Identification

When the goal of model selection is identification, the main drawback
of VFCV is that $n_t \ll n$ is often required for choosing
consistently the true model (see Section~\ref{sec.cvconsistent}),
whereas VFCV does not allow $n_t < n/2$.
Depending on the frameworks, different (empirical) recommandations
for choosing $V$ can be found in the literature. In ordered variable
selection, the largest $V$ seems to be the better, $V=10$ providing
results close to the optimal ones \citep{Zha:1993}. On the contrary,
\cite{Die:1998} and \cite{Alp:1999} recommend $V=2$ for choosing the
best learning procedures among two candidates.

\subsection{Future research} \label{sec.conclu.futur}

Perhaps the most important direction for future research would be to
provide, in each specific framework, precise quantitative measures of
the variance of CV estimators of the risk, depending on $n_t$, the
number of splits, and how the splits are chosen.
Up to now, only a few precise results have been obtained in this
direction, for some specific CV methods in linear regression or
density estimation (see Section~\ref{sec.riskestim.var}).
Proving similar results in other frameworks and for more general CV
methods would greatly help to choose a CV method for any given model
selection problem.

More generally, most theoretical results are not precise enough to
make any distinction between the hold-out and CV methods having the
same training sample size $n_t$, because they are equivalent at first
order. Second order terms do matter for realistic values of $n$,
which shows the dramatic need for theory that takes into account the
variance of CV when comparing CV methods such as VFCV and RLT with
$n_t = n(V-1)/V$ but $B \neq V$.

%
%

%%%%%%%%%%%%%%%%%%%%%%%%%%%%%%%%%%%%%%%%%%%%%%%%%%%%%%%%%%%%%%%%%%%%%%%
%%%%%%%%%%%%%%%%%%%%%%%%%%%%%%%%%%%%%%%%%%%%%%%%%%%%%%%%%%%%%%%%%%%%%%%
%%%%%%%%%%%%%%%%%%%%%%%%%%%%%%%%%%%%%%%%%%%%%%%%%%%%%%%%%%%%%%%%%%%%%%%
%%%%%%%%%%%%%%%%%%%%%%%%%%%%%%%%%%%%%%%%%%%%%%%%%%%%%%%%%%%%%%%%%%%%%%%
%%%%%%%%%%%%%%%%%%%%%%%%%%%%%%%%%%%%%%%%%%%%%%%%%%%%%%%%%%%%%%%%%%%%%%%

\bibliographystyle{apalike} %acmtrans-ims} %abbrvnat} %acmtrans-ims} %plainnat}
\bibliography{surveyCV,bibliosyl}

\end{document}